\documentclass[10pt,twocolumn,letterpaper]{article}

\usepackage{iccv}
\usepackage{times}
\usepackage{epsfig}
\usepackage{epstopdf}
\usepackage{graphicx}
\usepackage{enumerate}
\usepackage{amsmath}
\usepackage{amssymb}
\usepackage[ruled,vlined]{algorithm2e}
\usepackage{bbm}
\usepackage{listings}
\usepackage{subcaption}
\usepackage{kantlipsum} 
\usepackage{mwe} 

\iccvfinalcopy 

\def\iccvPaperID{3001} 

\ificcvfinal\fi
 
\begin{document}

\title{Weighted Low Rank Approximation for Background Estimation Problems}

\author{Aritra Dutta\\
	University of Central Florida\\
	4000 Central Florida Blvd, Orlando, FL-32816\\
	{\tt\small d.aritra2010@knights.ucf.edu}
	\and
	Xin Li\\
	University of Central Florida\\
	4000 Central Florida Blvd, Orlando, FL-32816\\
	{\tt\small xin.li@ucf.edu}
}

\maketitle

\begin{abstract}
Classical principal component analysis~(PCA) is 
not robust to the presence of sparse outliers in the data.~The use of the $\ell_1$ norm in the Robust PCA~(RPCA) method successfully eliminates the weakness of PCA in separating the sparse outliers.~In this paper, by sticking a simple weight to the Frobenius norm, we propose a weighted low rank~(WLR) method to avoid the often computationally expensive algorithms relying on the $\ell_1$ norm. As a proof of concept, a background estimation model has been presented and compared with two  $\ell_1$ norm minimization algorithms.~We illustrate that as long as a simple weight matrix is inferred from the data, one can use the weighted Frobenius norm and achieve the same or better performance.
\end{abstract}

\vspace{-0.3in}
\section{Introduction}
\vspace{-0.1in}
In image processing,~rank-reduced signal processing, computer vision, and in many other
engineering applications the classical principal component analysis~(PCA) is a successful tool~\cite{pca}.~However, it might lead to a degraded construction in some cases as it is not able to preserve any structure of the data matrix.~In 1987, Golub \etal~\cite{golub} were the first to consider a {\it constrained} low rank approximation problem of matrices to address this fundamental flaw in PCA:~Given $A=(A_1\;A_2)\in\mathbb{R}^{m\times n}$ with $A_1\in\mathbb{R}^{m\times k}$ and $A_2\in\mathbb{R}^{m\times (n-k)}$, find $A_G=(\tilde{B}_1\;\tilde{B}_2)$ such that
\vspace{-0.1in}
\begin{eqnarray}
(\tilde{B}_1\;\tilde{B}_2)=\arg\min_{\substack{B=(B_1\;B_2)\\B_1=A_1\\{\rm rank}(B)\le r}}\|A-B\|_F^2.\label{golub's problem}
\end{eqnarray}
~\\[-0.15in]
That is, Golub \etal required that the first few columns, $A_1,$ of $A$ must be preserved when one looks for a low rank approximation of $(A_1\;A_2).$~As in the standard low rank approximation~(which is equivalent to PCA), the constrained low-rank approximation problem of Golub \etal has a closed form solution.

Inspired by (\ref{golub's problem}) above and motivated by applications in which $A_1$ may contain noise, it makes more sense if we require $\|A_1-B_1\|_F$ small instead of asking for $B_1=A_1$. This leads us to consider the following problem: Let $\eta>0$, find $(\hat{B}_1\;\;\hat{B}_2)$ such that
\vspace{-0.12in}
\begin{eqnarray}\label{closeness problem}
(\hat{B}_1\;\;\hat{B}_2)=\arg\min_{\substack{B=(B_1\;B_2)\\\|A_1-B_1\|_F\leq \eta\\{\rm rank}(B)\le r}}\| A-B\|_F^2.
\end{eqnarray}
~\\[-0.15in]
Or, for a large parameter $\lambda$, consider
\vspace{-0.04in}
\begin{eqnarray}\label{unconstraint closeness}
\min_{\substack{B=(B_1\;B_2)\\{\rm rank}(B)\le r}}\{ \lambda^2\|A_1-B_1\|_F^2
+\| A_2-{B}_2 \|_F^2\}.
\end{eqnarray}
~\\[-0.22in]
As it turns out, (\ref{unconstraint closeness}) can be viewed as a generalized total least squares problem~(GTLS) and can be solved in closed form as a special case of weighted low-rank approximation with a rank-one weight matrix by using a single SVD of the given matrix $(\lambda A_1\;\;A_2)$~\cite{markovosky1,markovosky}. Using the closed form solutions, one can verify that the solution to (\ref{golub's problem}) is the limit case of the solutions to (\ref{unconstraint closeness}) as $\lambda \to\infty$. Thus, (\ref{golub's problem}) can be viewed as a special case when $\lambda = \infty$. Note that, problem (\ref{unconstraint closeness}) can also be cast as a special case of structured low rank problems with element-wise weights~\cite{markovosky3,markovosky4}.~More specifically, we observe that (\ref{unconstraint closeness}) is contained in the following more general point-wise weighted low rank~(WLR) approximation problem~\cite{markovosky1,markovosky,srebro}:
\vspace{-0.1in}
\begin{eqnarray}\label{hadamard problem}
\min_{\substack{X=(X_1\;X_2)\\{\rm r}(X)\le r}}\|\left(A-X\right)\odot W\|_F^2,
\end{eqnarray}
~\\[-0.15in]
where $W\in\mathbb{R}^{m\times n}$ is a weight matrix and $\odot$ denotes the Hadamard product.

The idea of working with a weighted norm is very natural in solving many engineering problems. The weighted low rank approximation problem was studied first with $W$ being an indicator weight for dealing with the missing data case and then for more general weight in machine learning, collaborative filtering, 2-D filter design, and computer vision non-rigid shape and motion from image streams~\cite{srebro,srebromaxmatrix,Buchanan,manton,lupeiwang,shpak,wibergjapan,wiberg,kanade}. Working with a weighted norm can be challenging, as there is no closed form solution in general. 

~In the past decade, one of the most prevalent approaches used in background estimation is to treat it as a matrix decomposition problem~\cite{Bouwmans2016}.~Given a sequence of $n$ video frames with each frame converted into a vector ${\mathbf a}_i\in {\mathbb R}^m$, $i=1,2,...,n$,~the data matrix $A=({\mathbf a}_1, {\mathbf a}_2, ... , {\mathbf a}_n)\in {\mathbb R}^{m\times n}$ is the concatenation of all the frame vectors. As the background is not expected to change much throughout the frames when the camera motion is small, it is assumed to be low rank~\cite{oliver1999}.~At the same time, the foreground is usually sparse if its size is relatively small compared to the frame size~\cite{candeslimawright,APG,LinChenMa}. Therefore, it is natural to consider a matrix decomposition problem by decomposing $A$ as the sum of its background and foreground:
\vspace{-0.1in}
\begin{align*}
	A=B+F,
\end{align*}
~\\[-0.22in]
where $B,F\in {\mathbb R}^{m\times n}$ are the background and foreground matrices, respectively.~Using the above idea,~in~\cite{LinChenMa,candeslimawright,APG}, the robust principal component analysis (RPCA) was introduced to solve the background estimation problem by considering the background frames, $B$, having a low-rank structure and the foreground $A-B$ being sparse:
\vspace{-0.08in}
\begin{equation}\label{rpca}
\min_B\{\|A-B\|_{\ell_1}+\lambda \|B\|_*\}.
\end{equation}
~\\[-0.2in]
But the RPCA model cannot take advantage of possible extra information on the background.~In~\cite{xin2015},~Xin \etal recently proposed a stronger model named as generalized fused Lasso~(GFL) for the situation where pure background frames are given as a supervised learning method. 
Assuming that some pure background frames are given and the data matrix $A$ can be written into $A=(A_1 ~A_2)$, where $A_1$ contains the given pure background frames, Xin \etal in~\cite{xin2015} proposed the following model of the unknown matrices $B$ and $F$: with $B=(B_1~B_2)$ and $F=(F_1 ~F_2)$ partitioned in the same way as in $A$, find $B$ and $F$ satisfying
\vspace{-0.15in}
\begin{align*}
	\min_{\substack{B,F\\B_1=A_1}}{\rm rank}(B)+\|F\|_{gfl},
\end{align*}
~\\[-0.15in]
where $\|\cdot\|_{gfl}$ denotes a norm that is a combination of $l_1$ norm and a local spatial total variation norm (to encourage connectivity of the foreground).~When $B_1\neq A_1$,~Xin \etal referred the model as unsupervised model. Indeed, \cite{xin2015} further simplified the above model by assuming ${\rm rank}(B)={\rm rank}(B_1)$.~Since $B_1=A_1$ and $A_1$ is given, so $r:={\rm rank}(B_1)$ is also given and thus, we can re-write the model of \cite{xin2015} as follows:
\vspace{-0.11in}
\begin{equation}\label{gfl}
	\min_{\substack{B = (B_1\;B_2)\\{\rm rank}(B)\le r\\B_1=A_1}} \|A-B\|_{gfl}.
\end{equation}
~\\[-0.15in]
It is obvious that, except in different norms, problem (\ref{gfl}) is a constrained low rank approximation problem as in (\ref{golub's problem}).~In this paper,~we propose an algorithm to solve~(\ref{hadamard problem}) as a standalone problem for a special family of weights $W=(W_1\;\mathbbm{1})$,~where $\mathbbm{1}$ is matrix of all ones.~As a proof of concept, we present a background estimation model using our WLR algorithm as it seems a natural fit to the problem.~In addition, we compare the performance of our proposed model with the RPCA and GFL algorithms in background estimation with static and dynamic background.~Our main focus in this paper is not to propose a background estimation model rather show how a properly weighted Frobenius norm can be made robust to the outliers similar to the $\ell_1$ norm.~For a comprehensive review of the most recent and traditional algorithms for solving background estimation problem, we refer the reader to~\cite{Bouwmans201431,Bouwmans2016,Sobral20144}.

{\bf Main Contributions}:~In this paper, we want to show that through a special weighted version of low-rank approximation problem~(\ref{hadamard problem}), and by learning the weight from the data, one can find a more robust and efficient approach to solve the background estimation problem as compare to the RPCA and GFL algorithms.~Our proposed model is as efficient as~\cite{xin2015}, but does not require any prior information~(see, for \eg Section 3.4).~More specifically we show:~(1) Instead of assuming the pure background frames are given, our model allows frames that are close to the background be used.~(2)~These approximate background frames are not given to us but learned from the data.~(3) Our experiments demonstrate that one might replace the computationally expensive $\ell_1$ norm as in RPCA and GFL algorithms by a weighted Frobenius norm and achieve a superior or at least comparable performance in detecting the foreground moving object.

\vspace{-0.1in}
\section{An Algorithm for WLR}
\vspace{-0.08in}
In this section, we propose an algorithm to solve~(\ref{hadamard problem}) for a special family of weights when $W=(W_1\;\mathbbm{1})$.~In fact even in this special case, our model shows superior performance in solving background estimation problem according to our experiments.~For convenience, let ${\rm r}(X_1)=k$.~Then any $X_2$ such that ${\rm r}(X_1\;\;X_2)\le r$ can be given in the form
\vspace{-0.1in}
\begin{align*}X_2=X_1C+BD,\end{align*}
~\\[-0.24in]
for some arbitrary matrices $B\in\mathbb{R}^{m\times (r-k)},$ $D\in\mathbb{R}^{(r-k)\times (n-k)},$ and $C\in\mathbb{R}^{k\times (n-k)}.$ Therefore,  problem (\ref{hadamard problem}) with $W=(W_1\;\mathbbm{1})$ is further reduced to:
\vspace{-0.1in}
\begin{align}\label{main problem 2}
	\min_{X_1,C,B,D}\left(\|(A_1-X_1)\odot W_1\|_F^2+\|A_2-X_1C-BD\|_F^2\right).
\end{align}
~\\[-0.15in]
Note that, for the special choice of the weight matrix, with a block structure $(X_1\;\;B)\begin{pmatrix} I_k & C\\0 & D\end{pmatrix},$ the problem~(\ref{main problem 2}) can be written alternatively 
in the framework of alternating weighted least squares algorithm in~\cite{markovosky}. Here we directly solve (\ref{main problem 2}) using a fast and simple numerical procedure based on the alternating direction method.  

Denote $F(X_1,C,B, D)=\|(A_1-X_1)\odot W_1\|_F^2+\|A_2-X_1C-BD\|_F^2$ as the objective function. The above problem (\ref{main problem 2}) can be numerically solved by using an alternating strategy~\cite{LinChenMa} of minimizing the function with respect to each component iteratively:
~\\[-0.23in]
\begin{eqnarray*}\label{update rule}
	\left\{\begin{array}{ll}
		\displaystyle{(X_1)_{p+1}=\arg\min_{X_1}F(X_1,C_p,B_p,D_p)},\\
		\displaystyle{C_{p+1}=\arg\min_{C}F((X_1)_{p+1},C,B_p,D_p)},\\
		\displaystyle{B_{p+1}=\arg\min_{B}F((X_1)_{p+1},C_{p+1},B,D_p)},\\
		\text{and,}\; \displaystyle{D_{p+1}=\arg\min_{D}F((X_1)_{p+1},C_{p+1},B_{p+1},D)}.
	\end{array}\right.
\end{eqnarray*}
\vspace{-0.02in}
\begin{algorithm}
	\SetAlgoLined
	\SetKwInOut{Input}{Input}
	\SetKwInOut{Output}{Output}
	\SetKwInOut{Init}{Initialize}
	\nl\Input{$A=(A_1\;\;A_2) \in\mathbb{R}^{m\times n}$ (the given matrix); $W= (W_1\;\;\mathbbm{1})\in\mathbb{R}^{m\times n}$ (the weight), threshold $\epsilon>0$\;}
	\nl\Init {$(X_1)_0,C_0,B_0,D_0$\;}
	\BlankLine
	\nl \While{not converged}
	{
		\nl $E_p=A_1\odot W_1\odot W_1+(A_2-B_pD_p)C_p^T$\;
		\BlankLine
		\nl \For {$i=1:m$}
		{	
			\nl $(X_1(i,:))_{p+1}=(E(i,:))_p({\rm diag}(W_1^2(i,1)$ $W_1^2(i,2)\cdots W_1^2(i,k))+C_pC_p^T)^{-1}$\;
		}
		\BlankLine
		\nl $C_{p+1}=((X_1)_{p+1}^T(X_1)_{p+1})^{-1}((X_1)_{p+1}^TA_2-(X_1)_{p+1}^TB_pD_p)$\;
		\nl $B_{p+1}=(A_2D_p^T-(X_1)_{p+1}C_{p+1}D_p^T)(D_pD_p^T)^{-1}$\;
		\nl $D_{p+1}=(B_{p+1}^TB_{p+1})^{-1}(B_{p+1}^TA_2-B_{p+1}^T(X_1)_{p+1}C_{p+1})$\;
		\nl $p=p+1$\;
	}
	\BlankLine
	\nl \Output{$(X_1)_{p+1}, (X_1)_{p+1}C_{p+1}+B_{p+1}D_{p+1}.$}
	\caption{WLR Algorithm}
\end{algorithm}
\vspace{-0.02in}
Each sub-problem above can be solved explicitly as described in {\bf Algorithm 1}. 

\begin{figure}
	\begin{center}
		\includegraphics[width=0.75\linewidth]{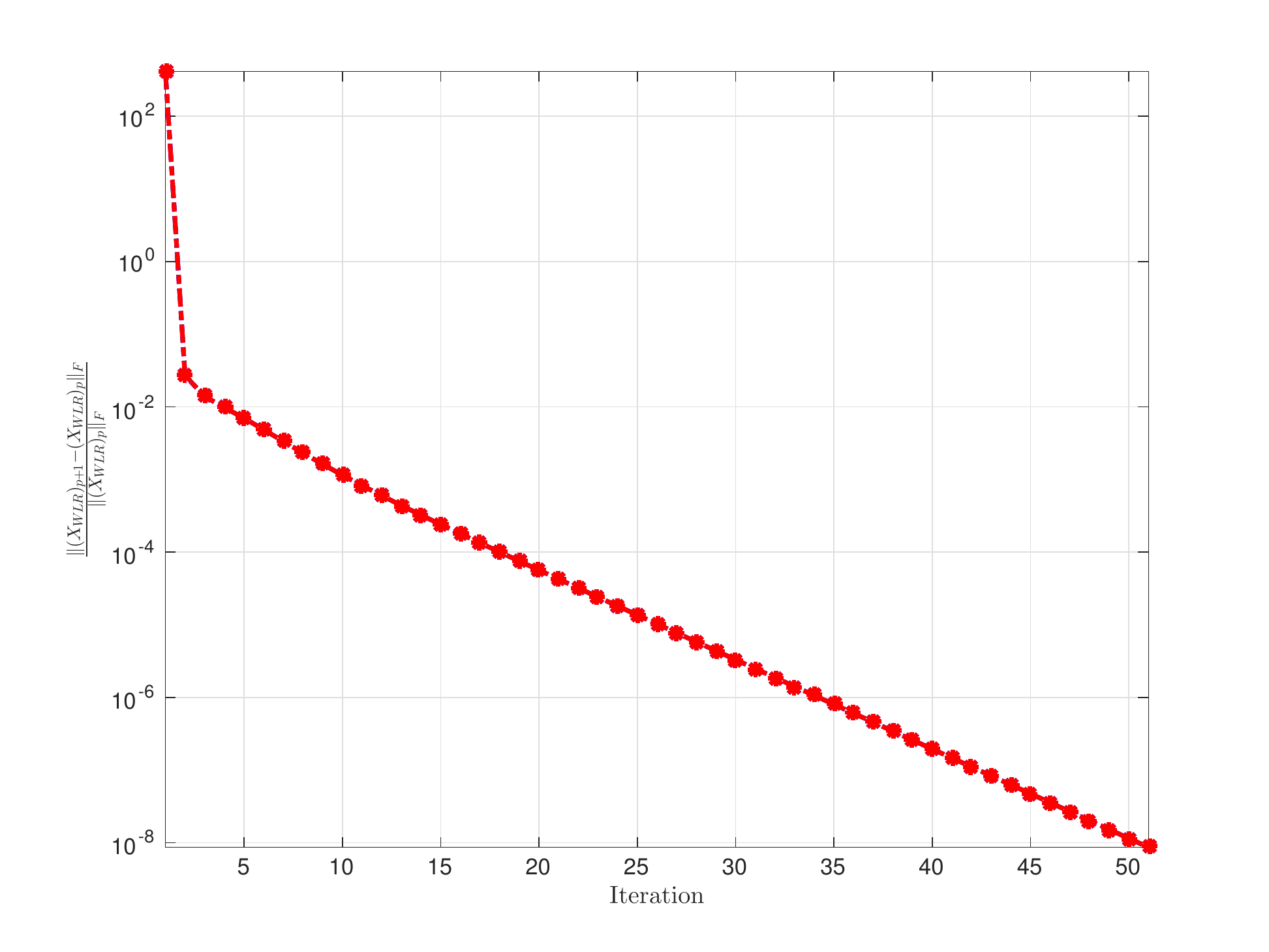}
	\end{center}
	\caption{Iterations vs. relative error on Stuttgart video sequence:~{\it Basic} scenario.}
	\label{conv}
\end{figure}

In our numerical procedure, we initialize $X_1$ and $D$ as random normal matrices and  $B$ and $C$ as zero matrices.~We denote $(X_{WLR})_p$ as our approximation to $A$ at $p$th iteration.~Using the notation we define $\|(X_{WLR})_{p+1}-(X_{WLR})_{p}\|_F=Error_p$ and as a measure of the relative error $\frac{Error_p}{\|(X_{WLR})_{p}\|_F}$ is used. For a threshold $\epsilon>0$ the stopping criteria of our algorithm at the $p$th iteration is  $Error_p<\epsilon$ or $\frac{Error_p}{\|(X_{WLR})_{p}\|_F}<\epsilon$ or if it reaches the maximum iteration.~Figure~\ref{conv} shows iteration $p$ vs. relative error plot for our algorithm on Stuttgart video sequence~(see Section 5 for more experimental detail) and it is clear from Figure~\ref{conv} that Algorithm 1 converges.~A detailed study of the convergence can be found in~\cite{duttali}.

\vspace{-0.15in}
\section{Background Estimation using WLR}
\vspace{-0.05in}
\begin{figure*}
	\centering
	\begin{subfigure}{.65\columnwidth}
		\includegraphics[width=\columnwidth]{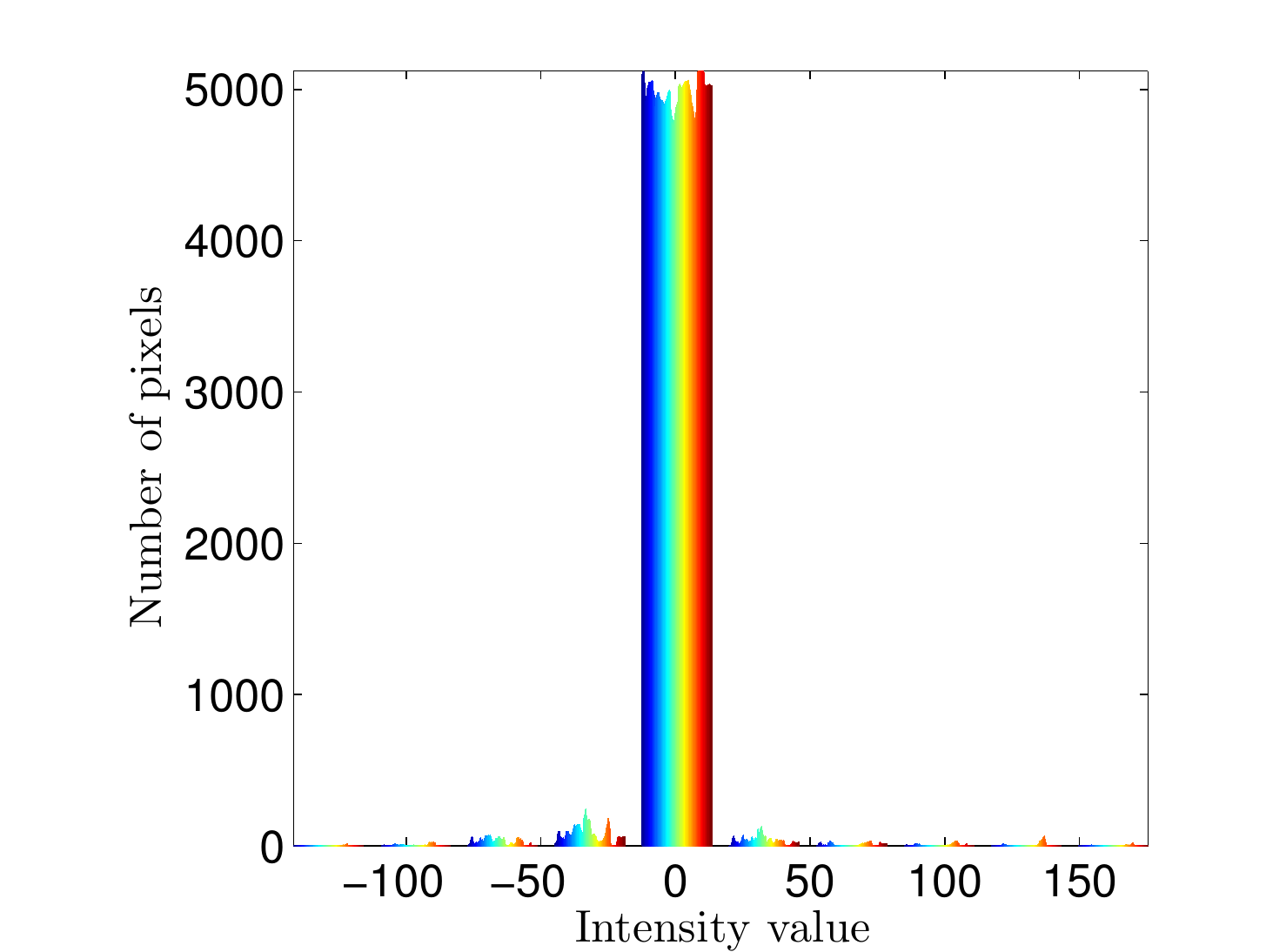}
		\caption{}
	\end{subfigure}
	\begin{subfigure}{.65\columnwidth}
		\includegraphics[width=\columnwidth]{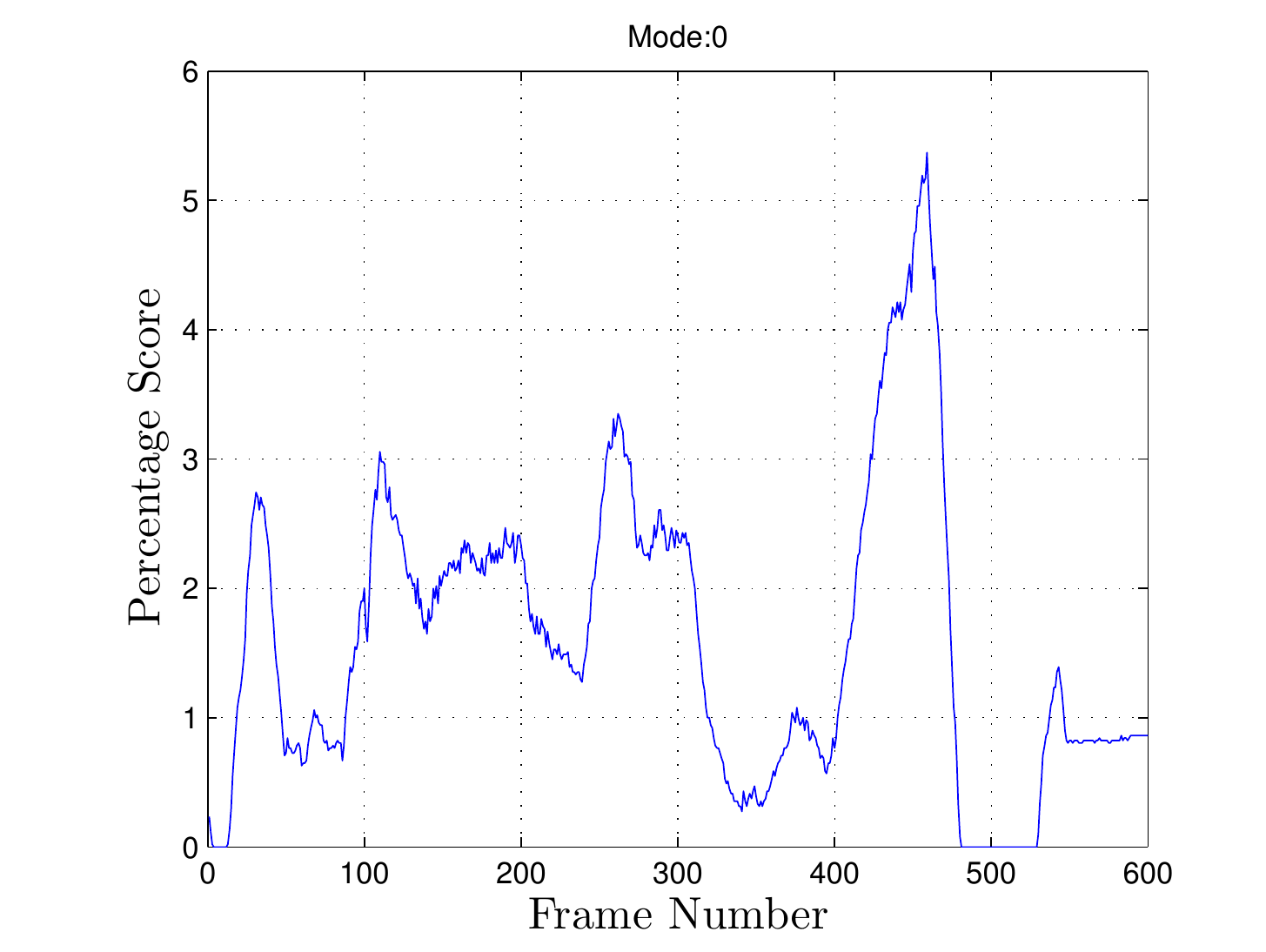}
		\caption{}
	\end{subfigure}
	\begin{subfigure}{.65\columnwidth}
		\includegraphics[width=\columnwidth]{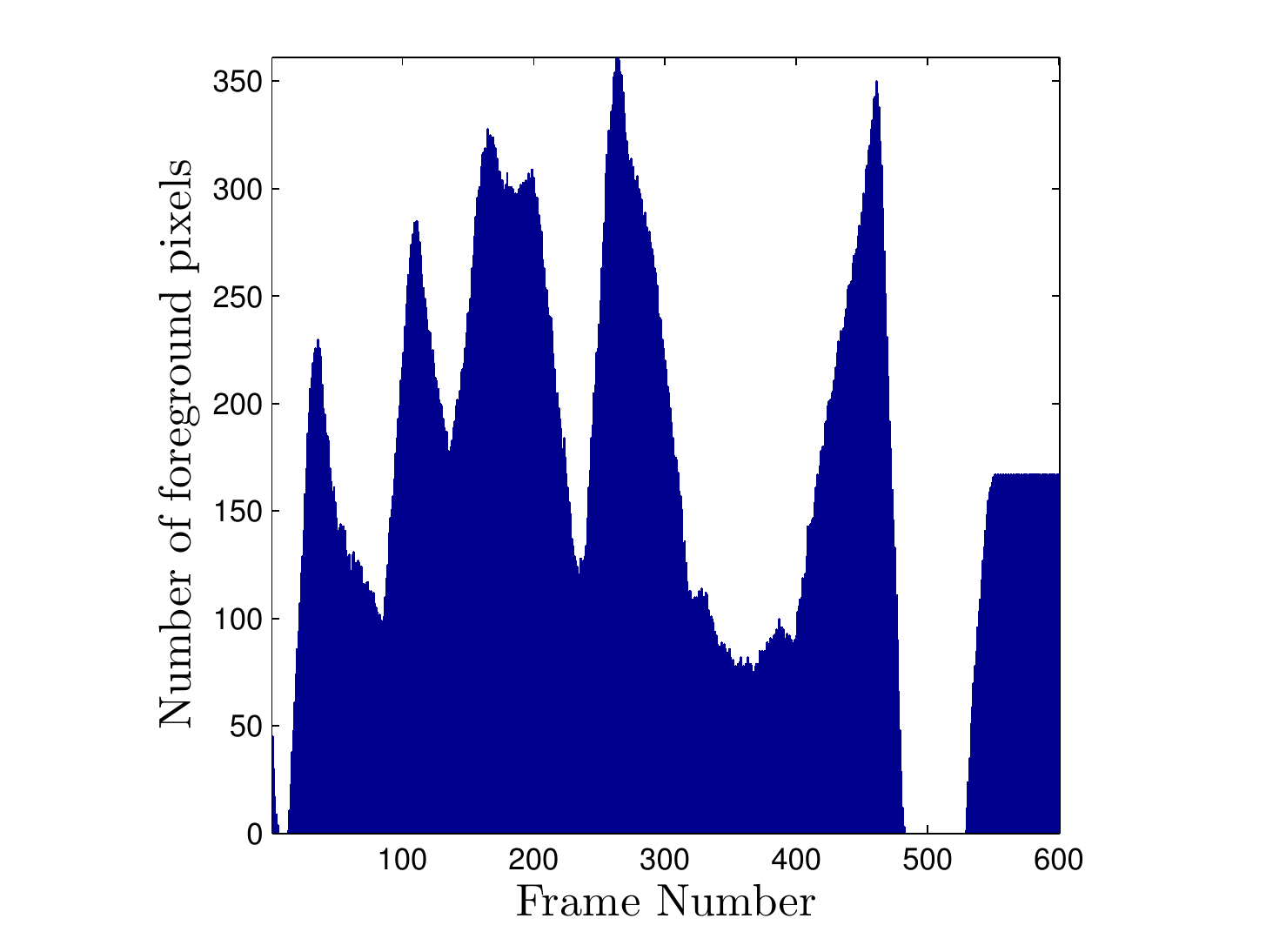}
		\caption{}
	\end{subfigure}
	\caption{Learning the weighted frame indexes for the {\it Basic} scenario using~\cite{duttaligongshah}.~(a)~Histogram to choose the threshold $\epsilon_1$.~(b)~Percentage score plot for 600 frames.~(c)~Original logical $G$ column sum, which indicates we are able pick up the indexes correctly corresponding to the frames that have least foreground movement. Originally, there are 53 frames in $G$ that have less than 5 pixels.~Using~\cite{duttaligongshah} we picked up 58 frame indexes on the~{\it Basic} scenario.}\label{learn_wt}
\end{figure*}

In this section,~we propose a background estimation model using Algorithm 1 and~show the power of our model over the existing RPCA and GFL algorithms.~To implement our proposed algorithm in the background estimation model, we use the heuristic to divide the data matrix $A$ into two blocks: $A_1$ and $A_2$, where $A_1$ mainly contains the background information, while $A_2$ contains both the background and foreground information.~We want to find a low-rank matrix $X=(X_1\;X_2)$ with compatible block partition, such that~$X_1\approx A_1$.~The pointwise multiplication with the weight matrix $W=(W_1\;~\mathbbm{1})$ helps us in that regard as $W_1\to\infty$. Finally, we point out that Xin \etal~ and Dutta \etal~\cite{duttali_acl,xin2015} also used the background frames in designing the weight matrix but they assumed that these background frame indexes are given. The novelty of our work is instead of using the prior knowledge of the available background frames, we {\it learn} the background frame indexes from the data and propose a robust background estimation model which is more realistic and applicable to real world problems.   
~\\[-0.2in]
\begin{algorithm}
	\SetAlgoLined
	\SetKwInOut{Input}{Input}
	\SetKwInOut{Output}{Output}
	\SetKwInOut{Init}{Initialize}
	\nl\Input{$A=(A_1\;\;A_2) \in\mathbb{R}^{m\times n}$ (the given matrix); $W= (W_1\;\;W_2)\in\mathbb{R}^{m\times n}, W_2=\mathbbm{1}\in\mathbb{R}^{m\times (n-k)}$ (the weight), threshold $\epsilon>0,$ $i_1,i_2\in\mathbb{N}$\;}
	
	\nl Run WSVT with $W=I_n$ to obtain: $A=B_{In}+F_{In}$\;
	\nl Plot image histogram of $F_{In}$ and find threshold $\epsilon_1$\;
	\nl Set $F_{In}(F_{In}\le\epsilon_1)=0$ and $F_{In}(F_{In}>\epsilon_1)=1$ to obtain a logical matrix $LF_{In}$\;
	\nl Convert $B_{In}$ directly to a logical matrix $LB_{In}$\;
	\nl Find
	$\epsilon_2={\rm mode}(\{\frac{\sum_i(LF_{IN})_{i1}}{\sum_i(LB_{IN})_{i1}}\times 100,\frac{\sum_i(LF_{IN})_{i2}}{\sum_i(LB_{IN})_{i2}}\times 100,\cdots, \frac{\sum_i(LF_{IN})_{in}}{\sum_i(LB_{IN})_{in}}\times 100\})
	$\;
	\nl Denote $S = \{i:(\frac{\sum_i(LF_{IN})_{i1}}{\sum_i(LB_{IN})_{i1}}\times 100,\frac{\sum_i(LF_{IN})_{i2}}{\sum_i(LB_{IN})_{i2}}\times 100,\cdots, \frac{\sum_i(LF_{IN})_{in}}{\sum_i(LB_{IN})_{in}}\times 100)\le \epsilon_2\}$\;
	\nl Set $k=\Bigl\lceil|S|/i_1\Bigr\rceil, r = k+i_2$\;
	\nl Rearrange data:~$\tilde{A}_1= (A(:,i))_{m\times k}$, $i\in S$ randomly chosen and $\tilde{A}_2= (A(:,i'))_{m\times (n-k)}$, $i\neq i'$\;
	\nl Apply {\bf Algorithm 1} on $\tilde{A}=(\tilde{A}_1\;\tilde{A}_2)$ to obtain $\tilde{X}$\;
	\nl Rearrange the columns of $\tilde{X}$ similar to $A$ to find $X$\;
	\nl \Output{$X$.}
	\caption{Background Estimation using WLR}
\end{algorithm}

\begin{figure}
	\begin{center}
		\includegraphics[width=\linewidth]{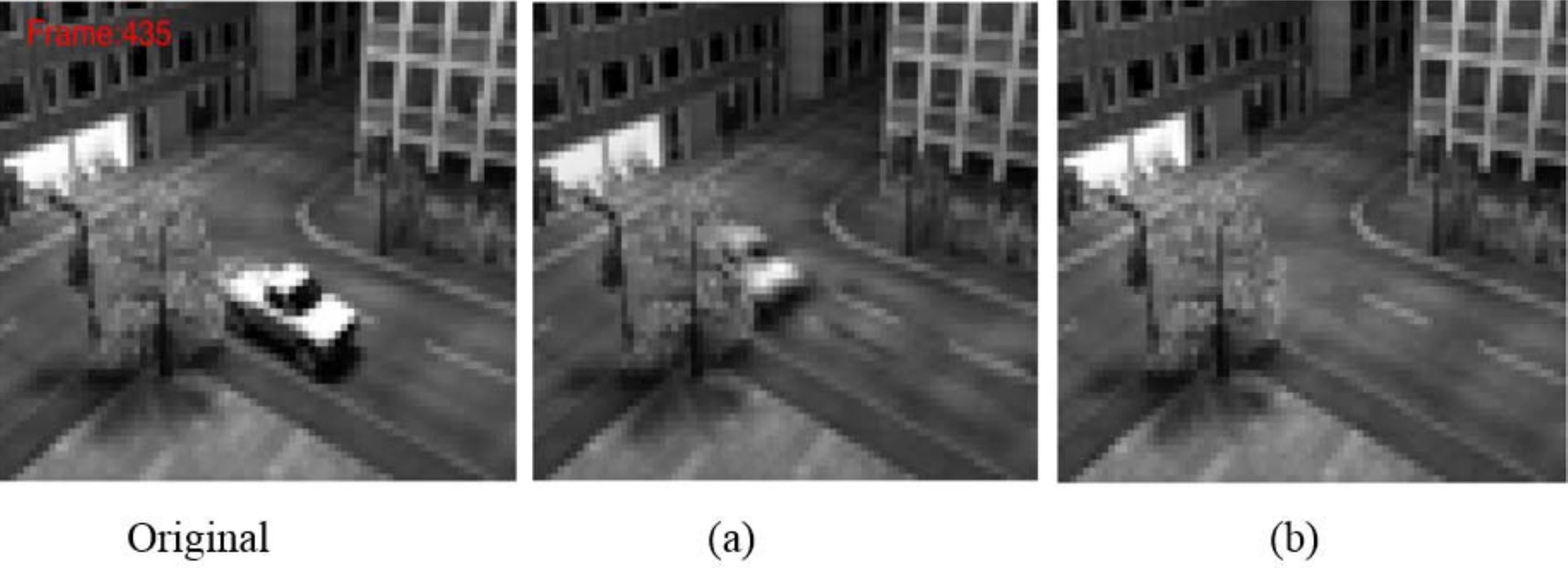}
	\end{center}
	\caption{The effect of using weights in WLR algorithm on the {\it Basic} scenario.~Frame number 435. Background estimation using WLR with:~(a) $(W_1)_{ij}\in[5,10]$,~(b)~$(W_1)_{ij}\in[500,1000].$~In (a) the estimated background has the remaining foreground object, but as we increase the weights, the foreground object disappears in (b).}
	\label{weight_435}
\end{figure}

\begin{table}
	\begin{center}
		\begin{tabular}{|l|c|c|c|}
			\hline
			Method    &{\it Basic} &{\it Noisy night}&{\it Light switch}\\
			\hline\hline
			WLR & {\bf 23.0676} & {\bf 24.0970} & {\bf 20.1874}\\
			iEALM & 160.251981 &108.679550 & 173.903928\\
			APG &107.982398 &115.547544 & 109.976457 \\
			\hline
		\end{tabular}
	\end{center}
	\caption{Average Computational time~(in seconds) for each algorithm in processing 600 frames of different scenarios.~All experiments were performed on a computer with 3.1 GHz Intel Core i7-4770S processor and 8GB memory.~The average computation time for iEALM and APG are almost 6.57 times and 4.95 times higher than that of WLR.}\label{table1}
\end{table}

In our experiments,~we extensively use~the Stuttgart synthetic video data set~\cite{cvpr11brutzer} for rigorous qualitative and quantitative comparisons.~It is a computer generated video sequence, that comprises both static and dynamic background/foreground objects and varying illumination in the background. 
We use three different test scenarios of the sequence:~(i) {\it Basic}: This scenario does not have noisy artifacts nor sudden illumination changes and is used as a general performance measure.~(ii) {\it Noisy night:} This scenario is a low-contrast nighttime scene, with increased sensor noise and camouflage.~(iii) {\it Light switch:} This scenario has varying illumination effects throughout the sequence.
Note that each scenario has 600 frames with identical foreground and background objects. Frame numbers 551 to 600 have static foreground, and frame numbers 6 to 12 and 483 to 528 have no foreground.~Additionally, the foreground comes with high quality ground truth mask available for each video frame.
\begin{figure}
	\centering  \includegraphics[width=\linewidth]{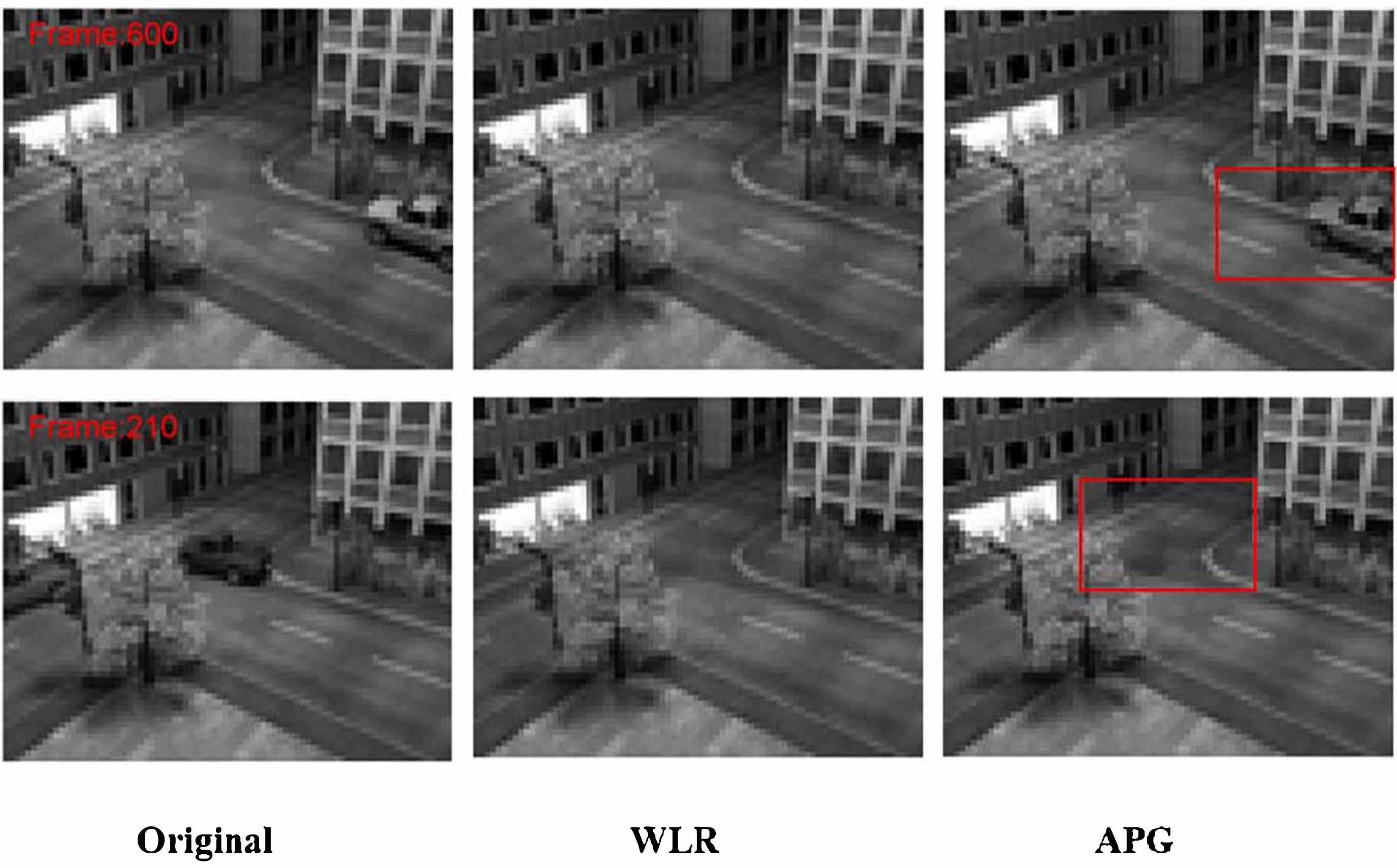}
	\caption{Background estimated by WLR and APG on the {\it Basic} scenario.~APG can not remove the static foreground object in frame 600.~On the other hand, in frame 210, the low-rank background estimated by APG has still some black patches. In both cases, WLR provides a substantially better background estimation than APG.}\label{basic}
\end{figure}
To compare with the existing RPCA algorithms, we use the inexact augmented Lagrange multiplier~(iEALM) method proposed by Lin \etal \cite{LinChenMa}, and the accelerated proximal gradient~(APG) algorithm proposed by Wright \etal \cite{APG}. For iEALM and APG, we set $\lambda={1}/{\sqrt{{\rm max}\{m,n\}}}$, and for iEALM we choose $\mu=1.5, \rho=1.25$ as used in \cite{LinChenMa,candeslimawright,APG}.
\begin{figure}
	\centering  \includegraphics[width=1.1\linewidth, height=2in]{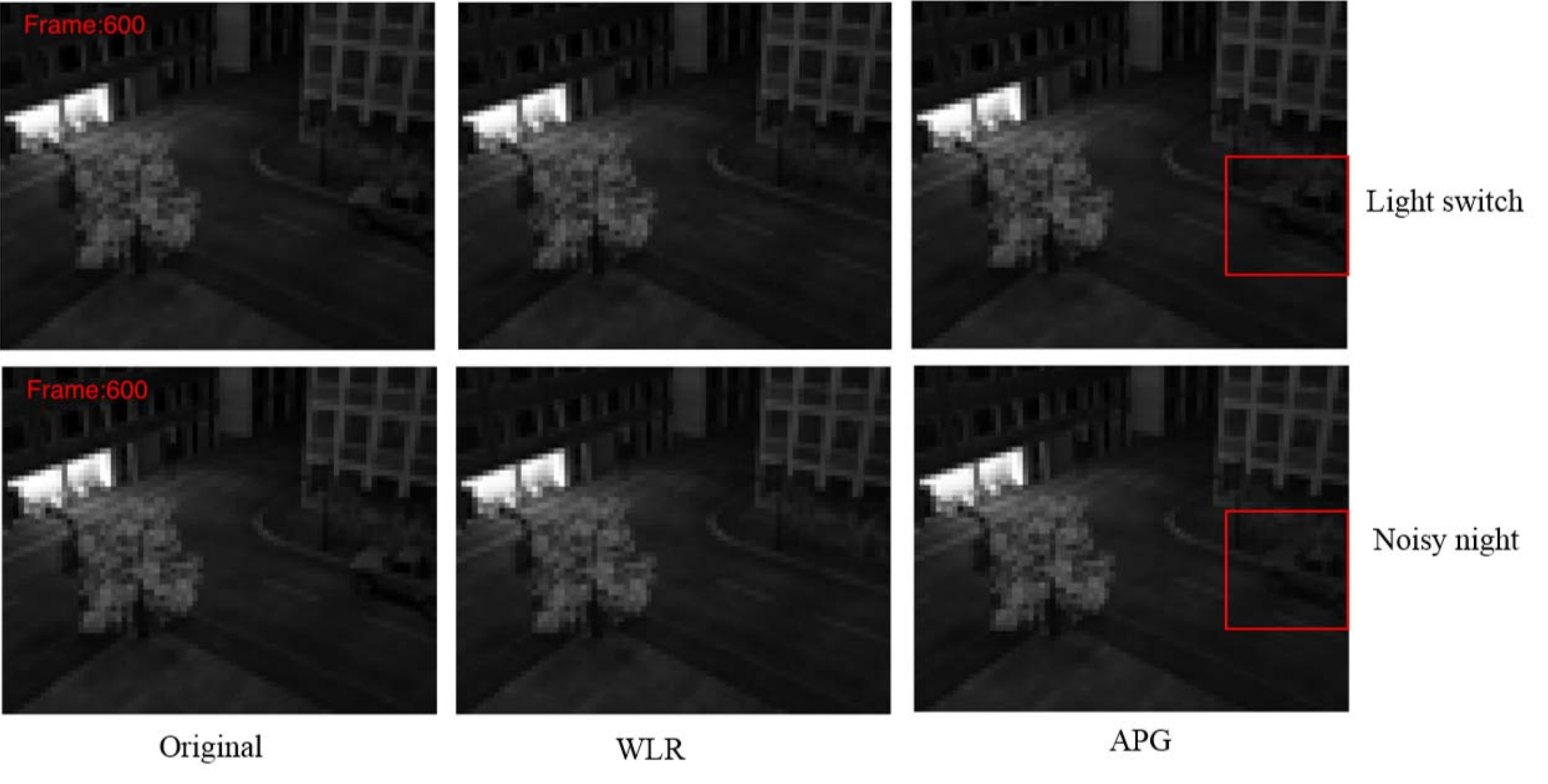}
	\caption{Background estimated by WLR and APG on {\it Light switch} and {\it Noisy night} scenario for frame 600. RPCA algorithm was not able to remove the static foreground object, but WLR provided an efficient background estimation by removing the static foreground object.}\label{nn_ls_600}
\end{figure}
~Given the sequence of 600 test frames, each frame in the test sequence is resized to $64\times80$; originally they were $600\times800$. Each resized frame is stacked as a column vector of size $5120\times1$ and we formed the test matrix $A$.~We denoted the ground truth matrix as $G$, with each column as a vectorized ground truth frame.~Then we apply the percentage score model described in~\cite{duttaligongshah} to learn the set $S$ that represents the frame indexes with least foreground movement based on a crude estimate of the initial background~($B_{In}$) and foreground~($F_{In}$) from the data matrix $A$. See Algorithm 2 for a detailed description of this method, and see Figure~\ref{learn_wt} for it's performance on the {\it Basic} scenario.~In our experiments, for the Stuttgart video sequence, we empirically choose $k=\Bigl\lceil|S|/2\Bigr\rceil$, where $|S|$ denotes the cardinality of the set $S$. We set $r=k+1$.~Therefore, following Algorithm 2, $i_1=2$ and $i_2=1$ for Stuttgart video sequence. 
\begin{figure}
	\centering  \includegraphics[width=\linewidth, height=1.25in]{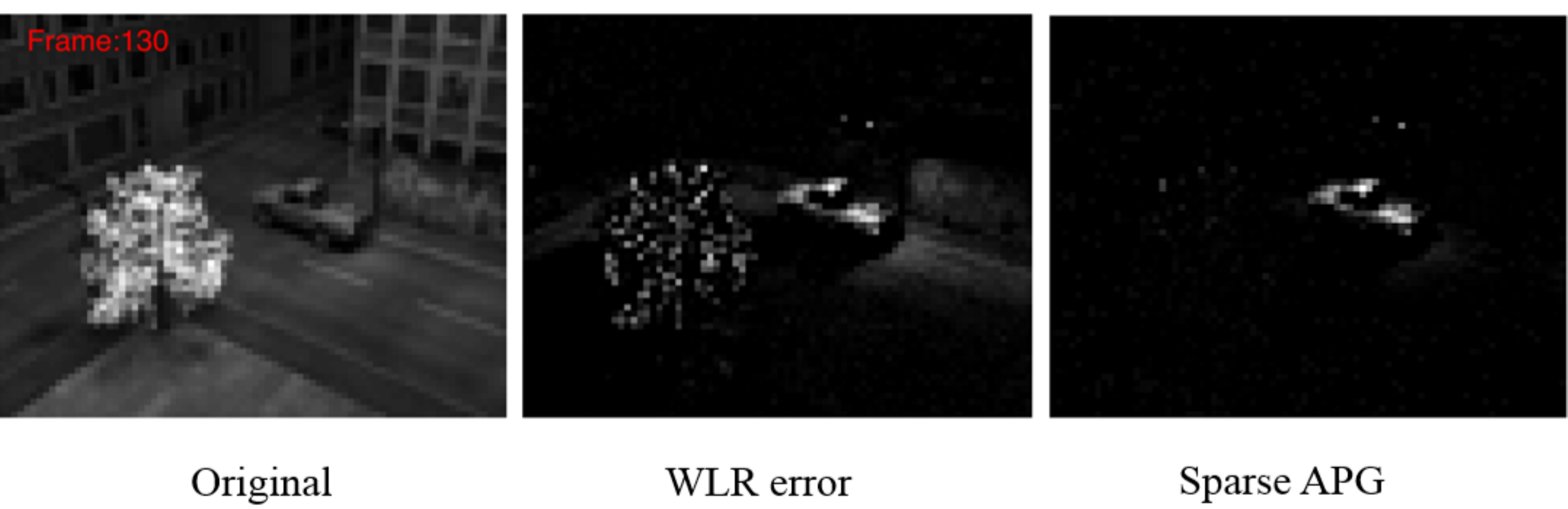}
	\caption{Foreground recovered by WLR and APG on the {\it Light switch} scenario, frame 130.~Starting from frame 125 the illumination changes suddenly. The sparse foreground recovered by APG does not capture the change in illumination.~WLR captures the effect of change in illumination, irregular movements of the tree leaves, and reflections effectively.}\label{ls_130_sp}
\end{figure}
\begin{figure}
	\centering  \includegraphics[width=1.1\linewidth, height=2in]{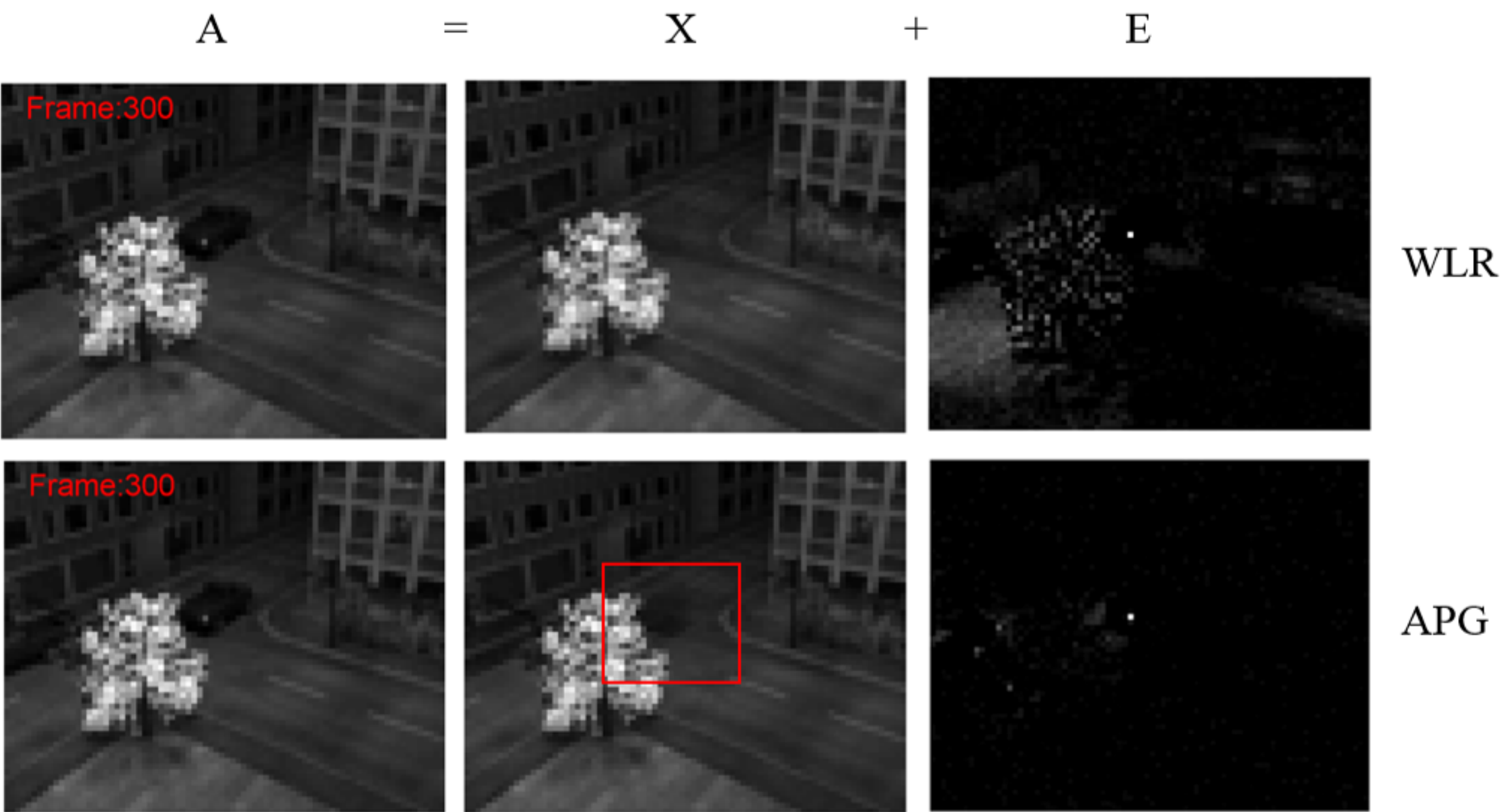}
	\caption{Background and foreground estimated by WLR and APG on {\it Light switch} scenario, frame 300. WLR has least MSSIM for frame 300 but still it provides a better visual quality foreground and background estimation than APG.~The red bounding box in APG frame is indicating the presence of foreground patch.}\label{nn_300}
\end{figure}
However, such assumptions do not apply to all practical scenarios.~Therefore, we argue that, in practical scenarios, the choices of $r$ and $k$ are problem-dependent and highly heuristic.~We rearrange the columns of our original test matrix $A$ as follows: Form $\tilde{A}_1= (A(:,i))_{m\times k}$ such that the indexes $i$ are randomly chosen from the set $S$, and form the second block $\tilde{A}_2$ using the remaining columns of the matrix $A$.~With the rearranged matrix $\tilde{A}=(\tilde{A}_1\;\;\tilde{A}_2)$ as our data matrix, we run Algorithm 1 for 50 iterations and obtain a low-rank estimation $\tilde{X}$. Finally, we rearrange the columns of $\tilde{X}$ as they were in the original matrix $A$ and form $X$. A threshold $\epsilon = 10^{-7}$ was chosen for Algorithm 1.

\begin{figure*}
	\centering
	\begin{subfigure}{.65\columnwidth}
		\includegraphics[width=\columnwidth]{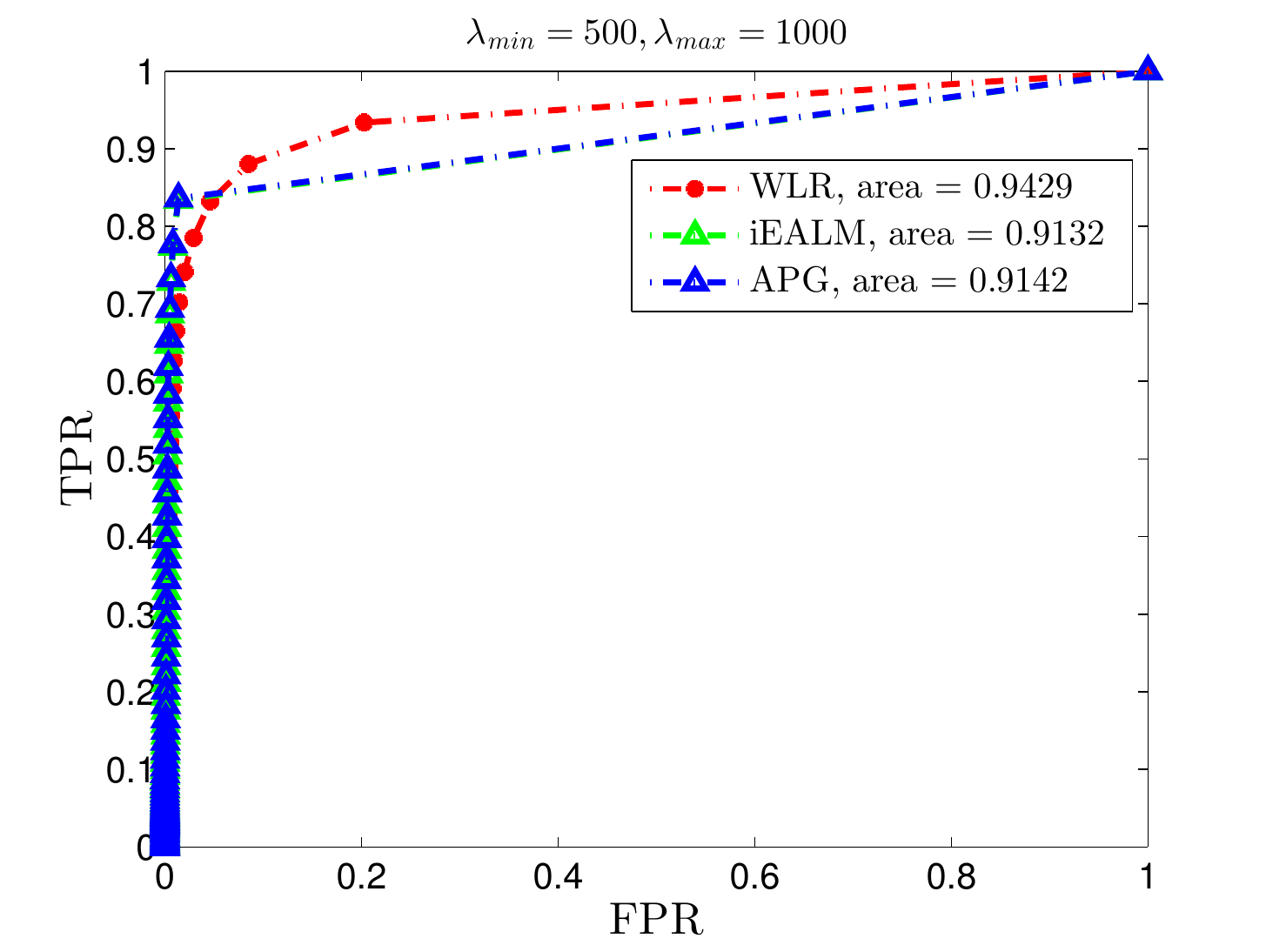}%
		\caption{Basic}%
	\end{subfigure}\hfill%
	\begin{subfigure}{.65\columnwidth}
		\includegraphics[width=\columnwidth]{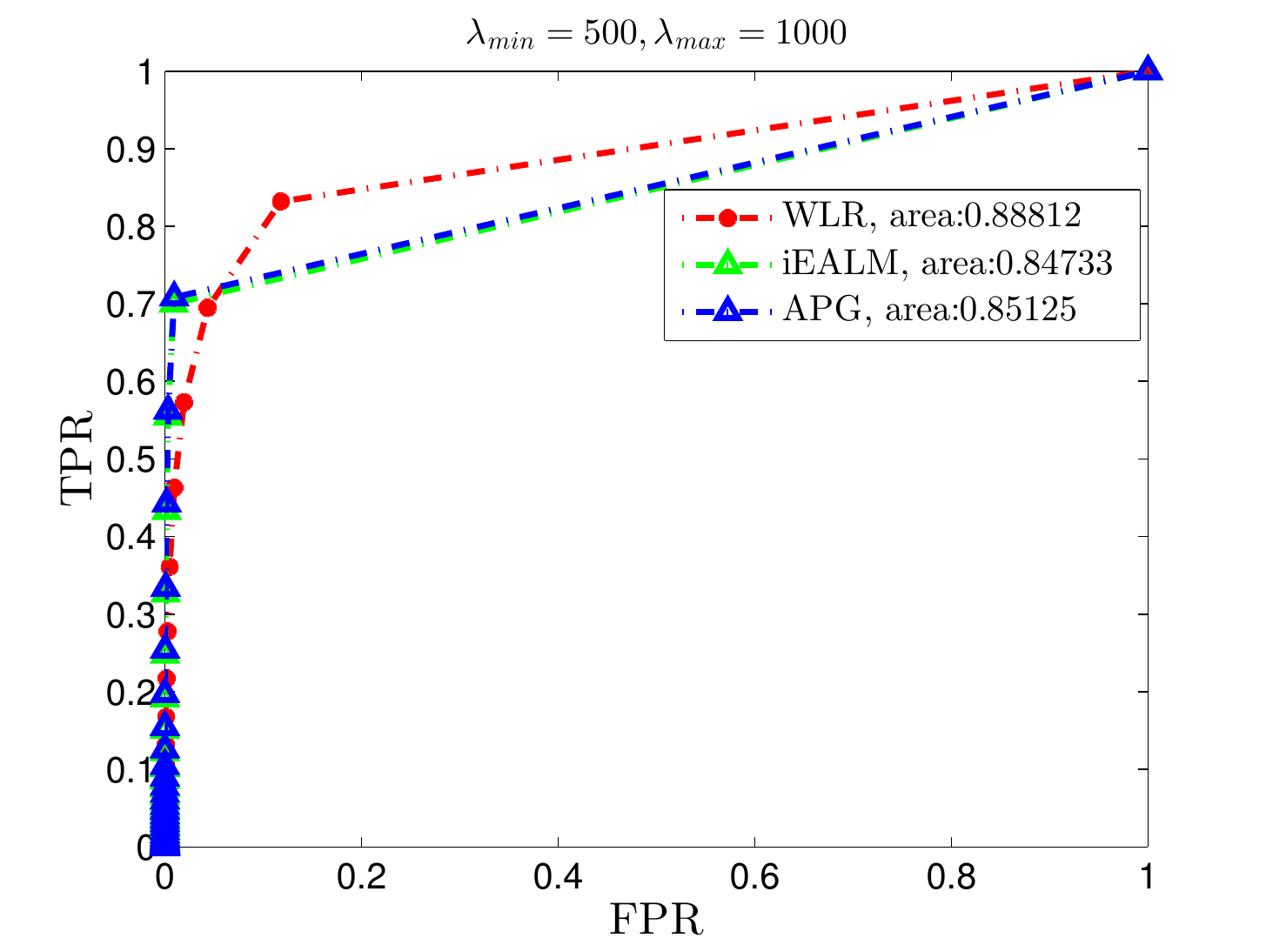}%
		\caption{Noisy night}%
	\end{subfigure}\hfill%
	\begin{subfigure}{.65\columnwidth}
		\includegraphics[width=\columnwidth]{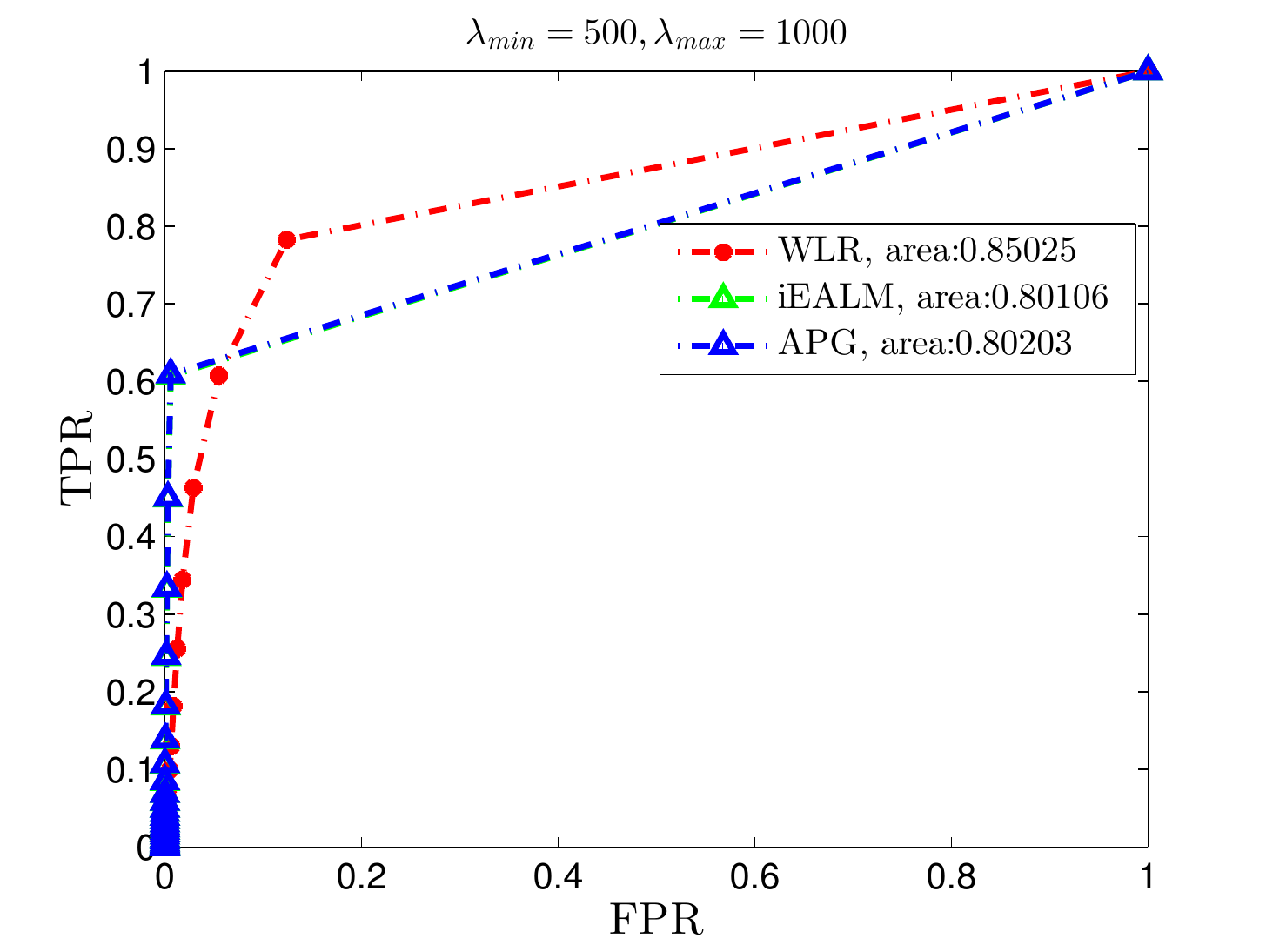}%
		\caption{Lightswitch}%
	\end{subfigure}%
	\caption{~ROC curve to compare between WLR,~iEALM, and APG.~The performance gain by WLR compare to APG on {\it Basic, Noisy night,} and {\it Light switch} scenarios are 3.252\%, 4.3313\%, and 6.012\% respectively, and compare to iEALM are 3.139\%, 4.8139\%, and 6.141\% respectively.}
	\label{roc_curve}
\end{figure*}

\vspace{-0.1in}
\subsection{Qualitative Analysis}
\vspace{-0.1in}
Since the background recovered by APG and iEALM have similar visual quality, we will only compare APG in this section.
\begin{figure}
	\begin{center}
		\includegraphics[width=0.7\linewidth]{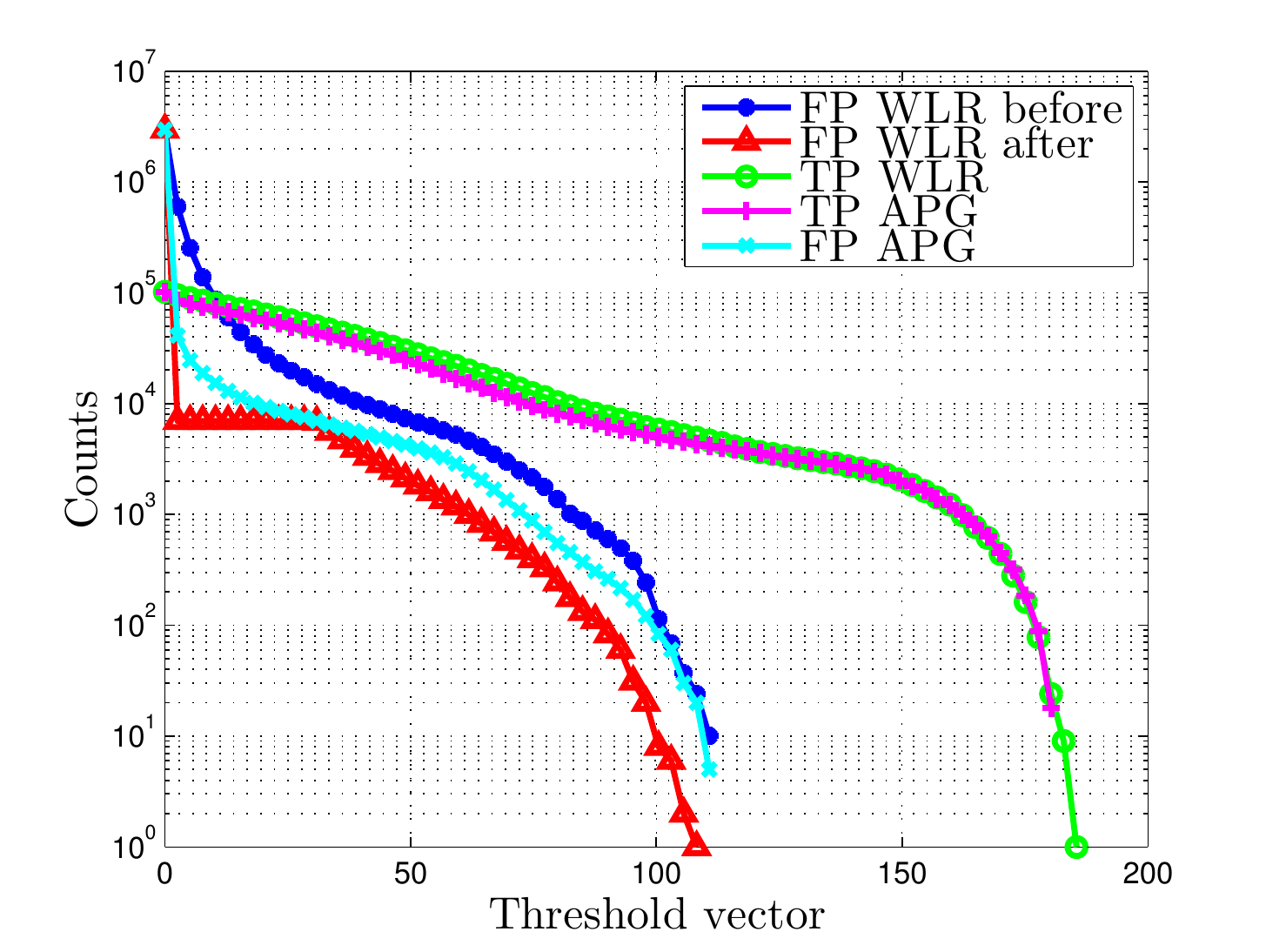}
	\end{center}
	\caption{True positive and false positive count for WLR and APG on~{\it Basic} scenario.~The false positive count for WLR substantially drops after thresholding $F$ by $\epsilon_1$.~On the other hand, WLR always has more or equal number of true positive count as APG.}
	\label{fp}
\end{figure}
\begin{figure}
	\centering
	\includegraphics[width=0.8\columnwidth]{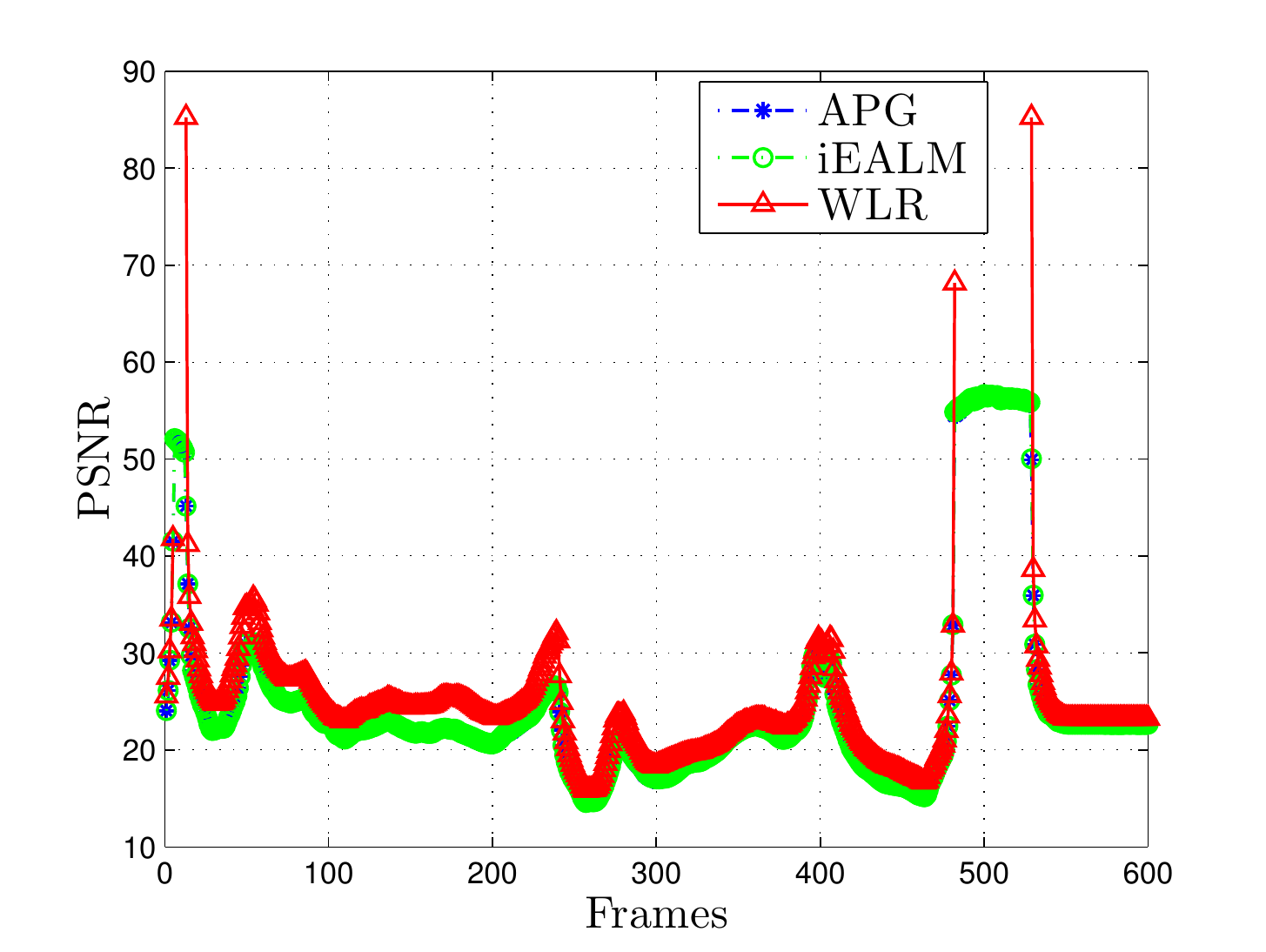}%
	\caption{Frames vs.~PSNR for {\it Basic} scenario.~The mean PSNR of APG and iEALM on the {\it Basic} sequence are 25.0092 and 25.0551,~respectively.~For WLR, the frames that do not contain the foreground object have PSNR equal to infinity, in all three scenarios.}
	\label{psnr}
\end{figure}


We present frame number 435 of the {\it Basic} scenario in~Figure~\ref{weight_435} to show the effect of a large weight,~$W_1$,~on the first block $A_1$:~our weighted low-rank algorithm can perform well in background estimation with proper choice of weight.~Next, in Figure~\ref{basic}, we present frame number 210 and 600 of the {\it Basic} scenario.~The performance of APG on frame 210 is comparable with WLR, but on frame 600, WLR clearly outperforms APG.~Finally, the experimental result in Figure~\ref{nn_ls_600}~shows the same phenomenon:~WLR completely removes the static foreground and provides a better visual background.~To conclude,~when the foreground is static, with the proper choice of $W, r,$ and $k$ our algorithm can provide a good estimation of the background by removing the static foreground object.
Our quantitative results in the next section suggest that RPCA algorithms act as a low-pass filter in the presence of a static foreground and attenuate its pixel values.~Therefore,~the pixels corresponding to the static foreground object stay as a part of the low-rank background. In presence of a ground truth mask of the foreground, the pixels of the frames corresponding to the static foreground captured as sparse components provide a comparable quantitative measure but create a poor human visual perception.~Figure~\ref{ls_130_sp} and~\ref{nn_300} present the foreground recovered by WLR and APG on the {\it Light switch} scenario and {\it Noisy night} scenario, respectively.
\begin{figure}
	\centering
	\includegraphics[width=0.8\columnwidth]{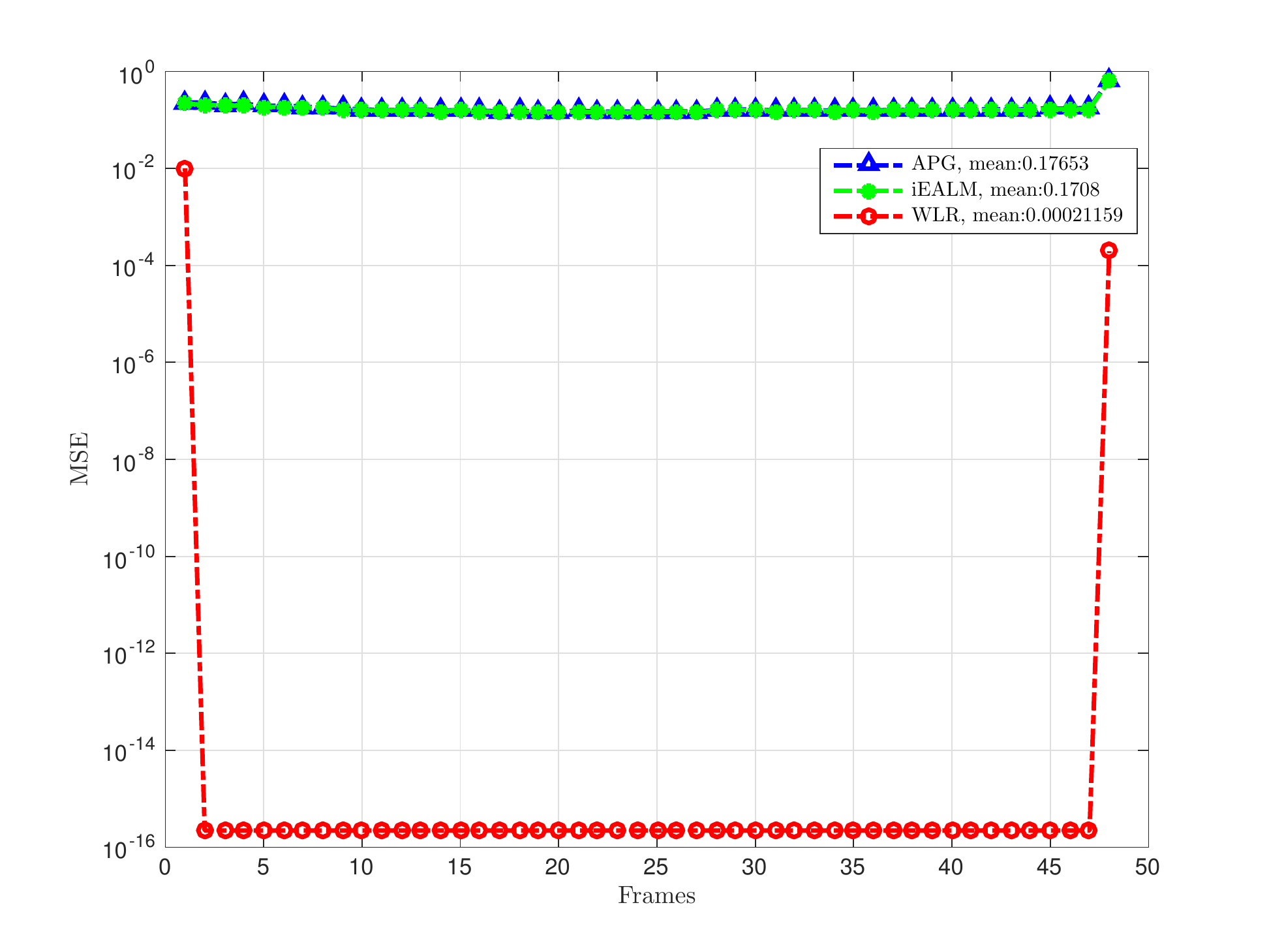}%
	\caption{Frames vs. MSE of different methods for frame numbers 482 to 529 for {\it Basic} scenario. MSE of WLR for frames 483 to 529  with no foreground movement is 0.}
	\label{mse}
\end{figure}
\begin{figure*}
	\centering
	\begin{subfigure}{.65\columnwidth}
		\includegraphics[width=\columnwidth]{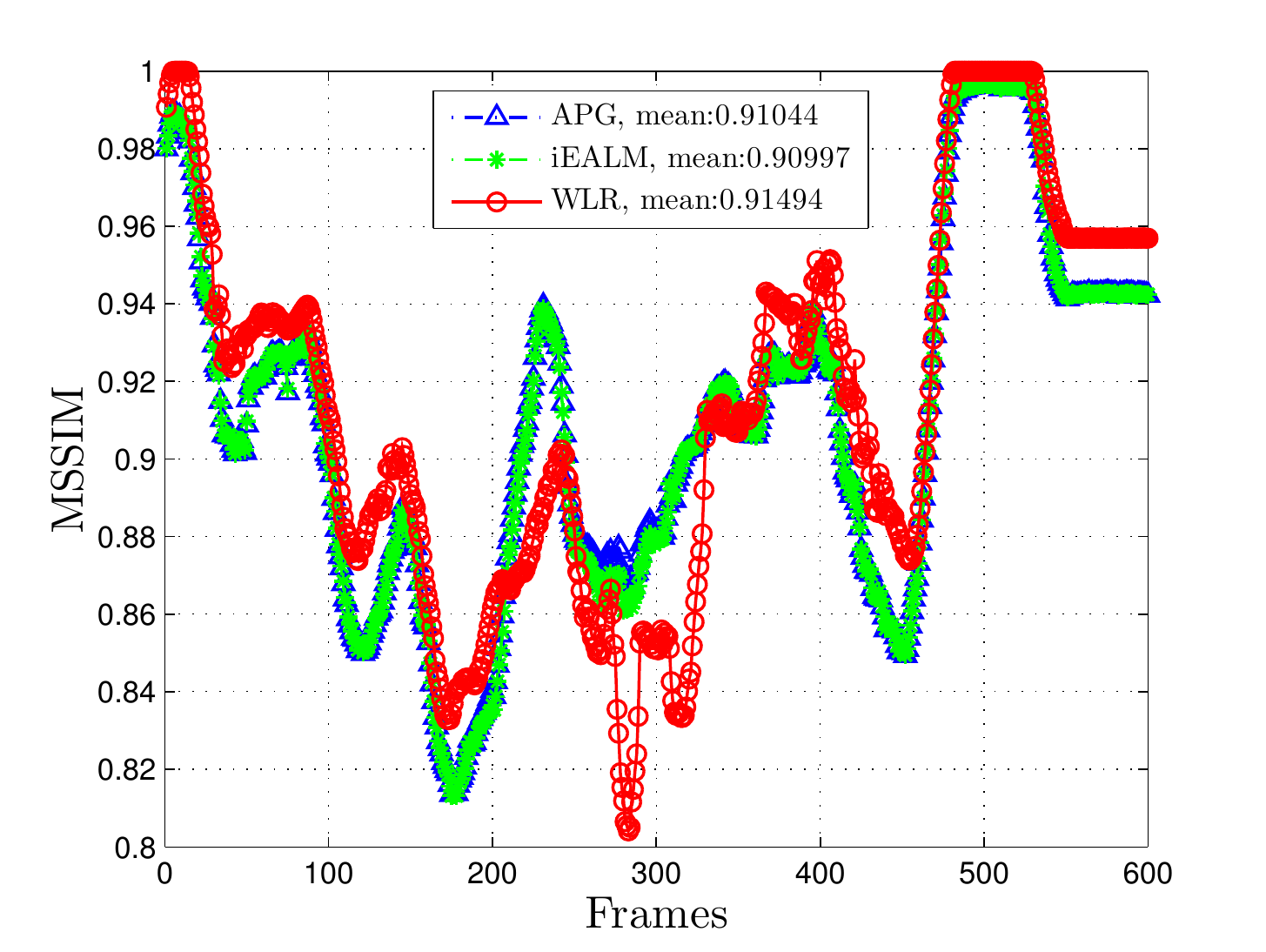}%
		\caption{Basic}%
	\end{subfigure}\hfill%
	\begin{subfigure}{.65\columnwidth}
		\includegraphics[width=\columnwidth]{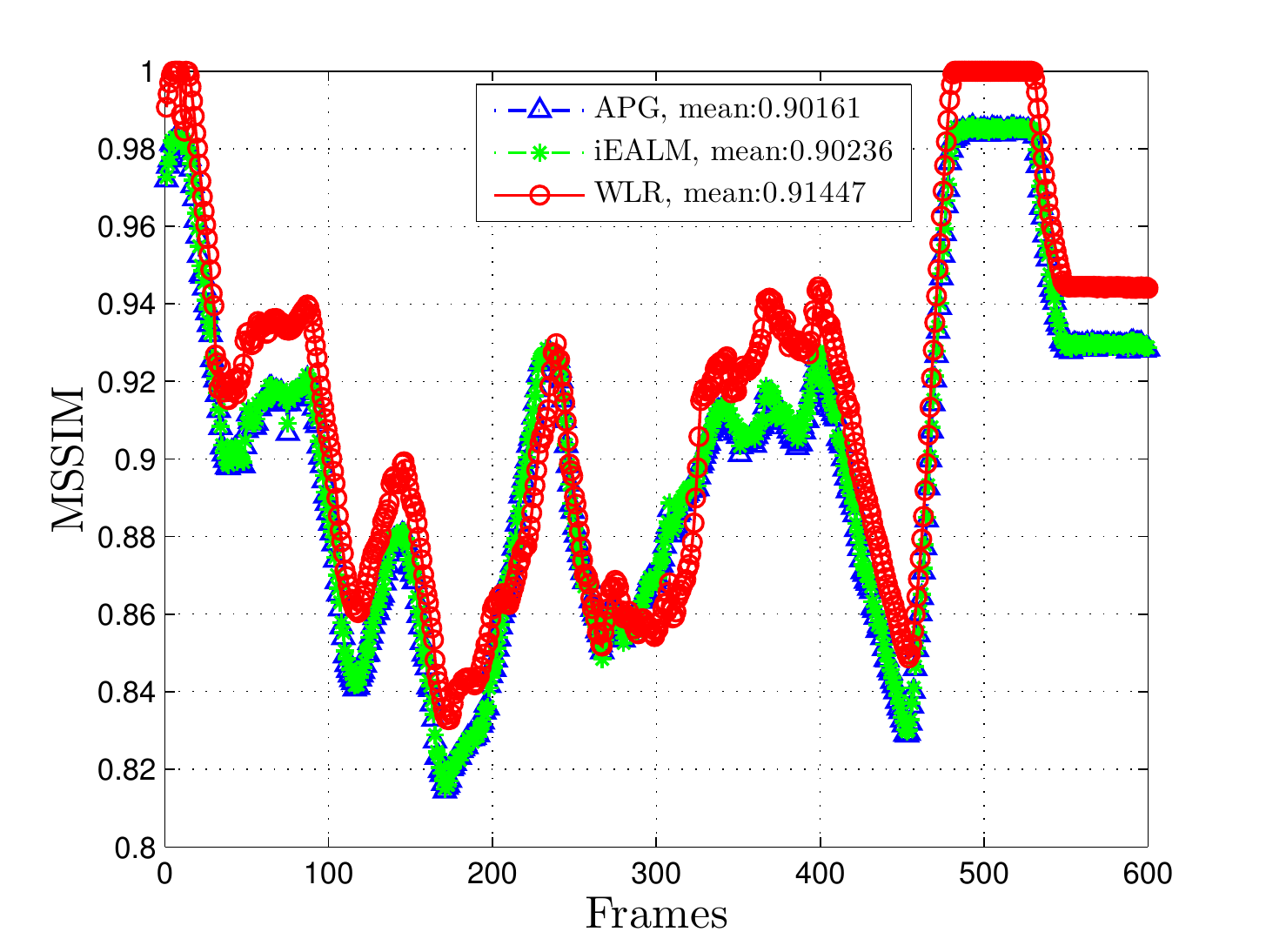}%
		\caption{Noisy night}%
	\end{subfigure}\hfill%
	\begin{subfigure}{.65\columnwidth}
		\includegraphics[width=\columnwidth]{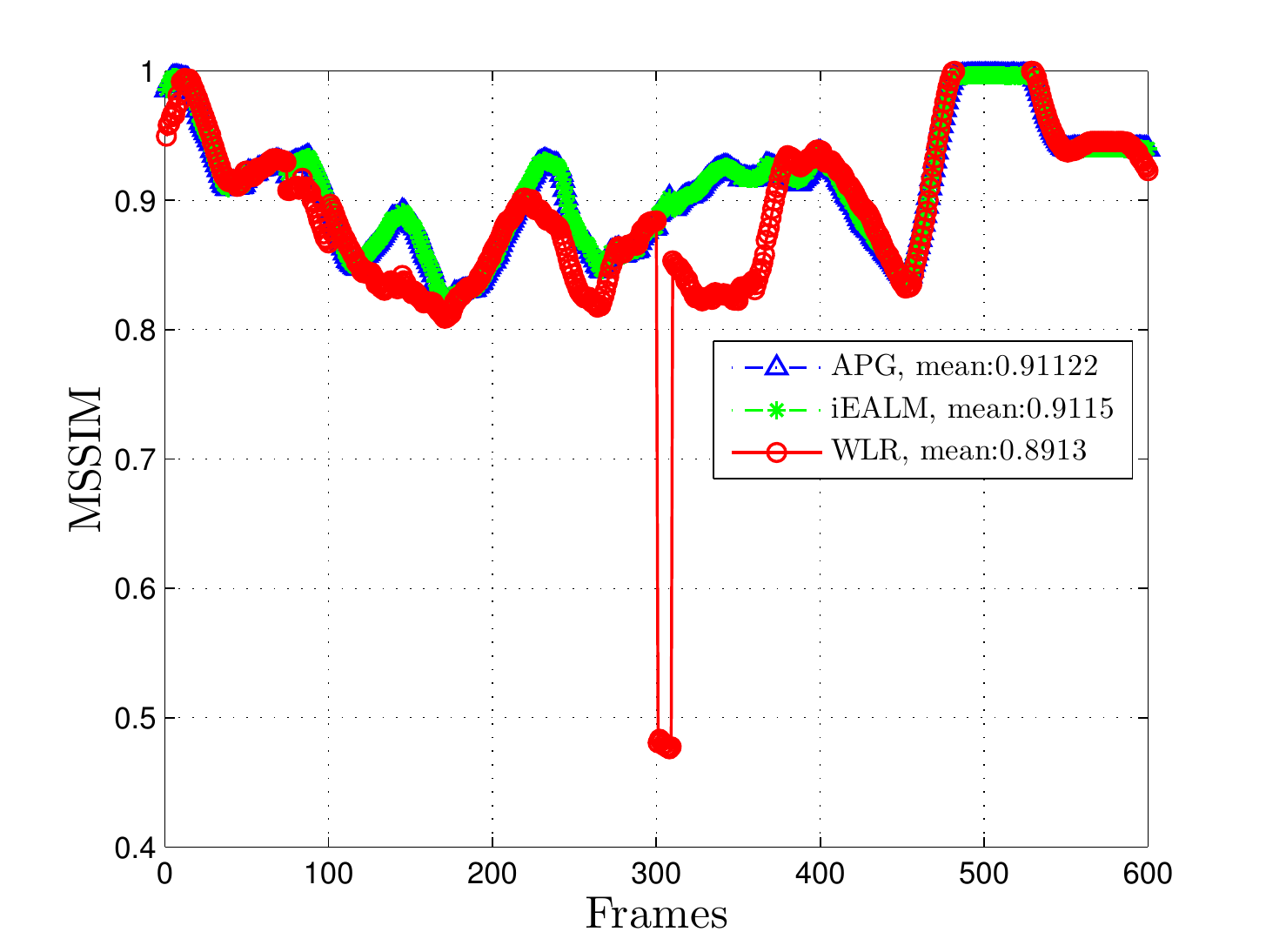}%
		\caption{Lightswitch}%
		\label{subfigc}%
	\end{subfigure}%
	\caption{MSSIM of different methods on all three scenarios.~WLR has better MSSIM compare to the RPCA algorithms corresponding to the frames which has static foreground or no foreground. The slight deterioration of performance of WLR in~(c) can be attributed by same choice of $i_1$ and $i_2$ in Algorithm 2 for all three scenarios.}
	\label{mssim}
\end{figure*}
We show WLR can capture the changing illumination and irregular dynamic background movements better than APG and can provide a visually better background frame, even on Frame number 300 of {\it Noisy night} scenario where WLR has least MSSIM. This can be attributed to the fact that,~RPCA algorithms are based on the assumption that the low-rank component is exactly low-rank while the sparse component being
exactly sparse~\cite{candeslimawright,APG,Bouwmans2016}.~The use of $\ell_1$ norm is good for removing the sparse components from the data but is not very capable of removing other artifacts.~Therefore, in the real-time video surveillance when the data are often corrupted by noisy artifacts, the assumption imposed on RPCA does not hold good. Considering the computational time of each algorithm from Table~\ref{table1}, WLR has minimal execution cost in producing a superior background estimation.
\vspace{-0.05in}
\subsection{Quantitative Analysis}
\vspace{-0.05in}
We now present different quantitative comparisons between the performance of our algorithm and that of the existing RPCA algorithms.~We use three different quantitative measures for this purpose: traditionally used receiver and operating characteristic~(ROC) curve, peak signal to noise ratio~(PSNR), and the most advanced measure mean structural similarity index~(MSSIM).
\begin{figure}
	\centering  \includegraphics[width=\linewidth, height=1.75in]{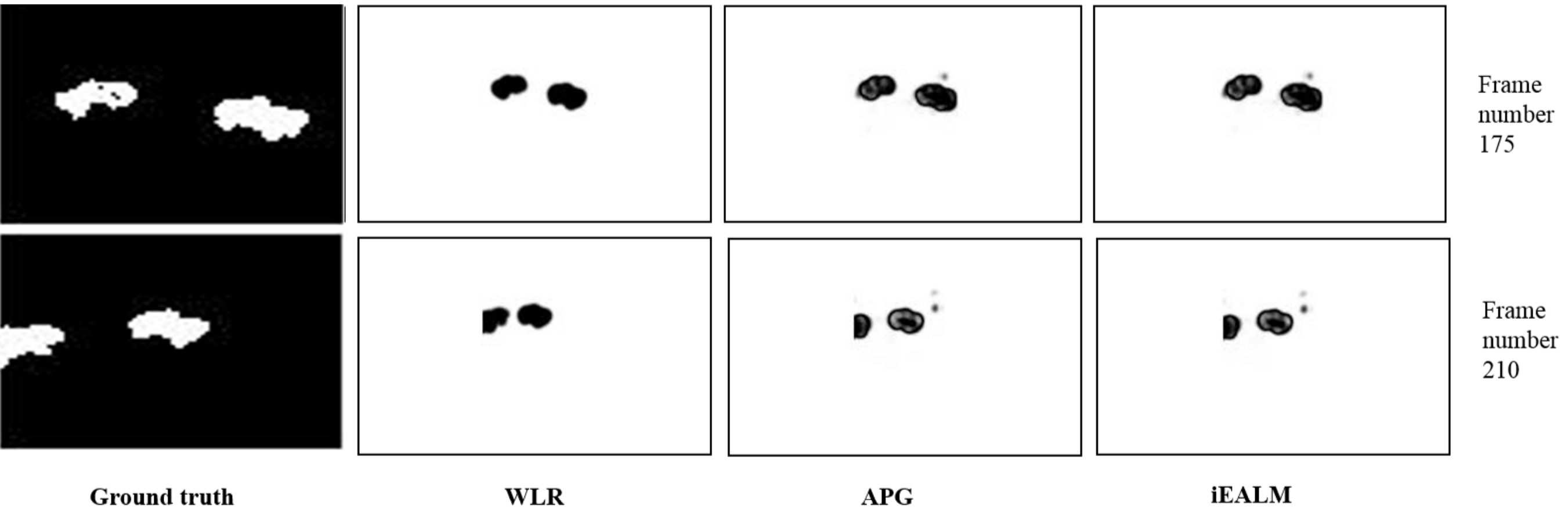}
	\caption{SSIM map for frame 175 and 210 of the {\it Basic} scenrio. Left to right: Ground truth frame~(size $64\times 80$), SSIM index map~(size $54\times 70$)~for WLR, APG, and iEALM.~WLR has superior SSIM index map than RPCA algorithms.}\label{ssim_map}
\end{figure}
\begin{figure}
	\centering  \includegraphics[width=\linewidth, height=2.9in]{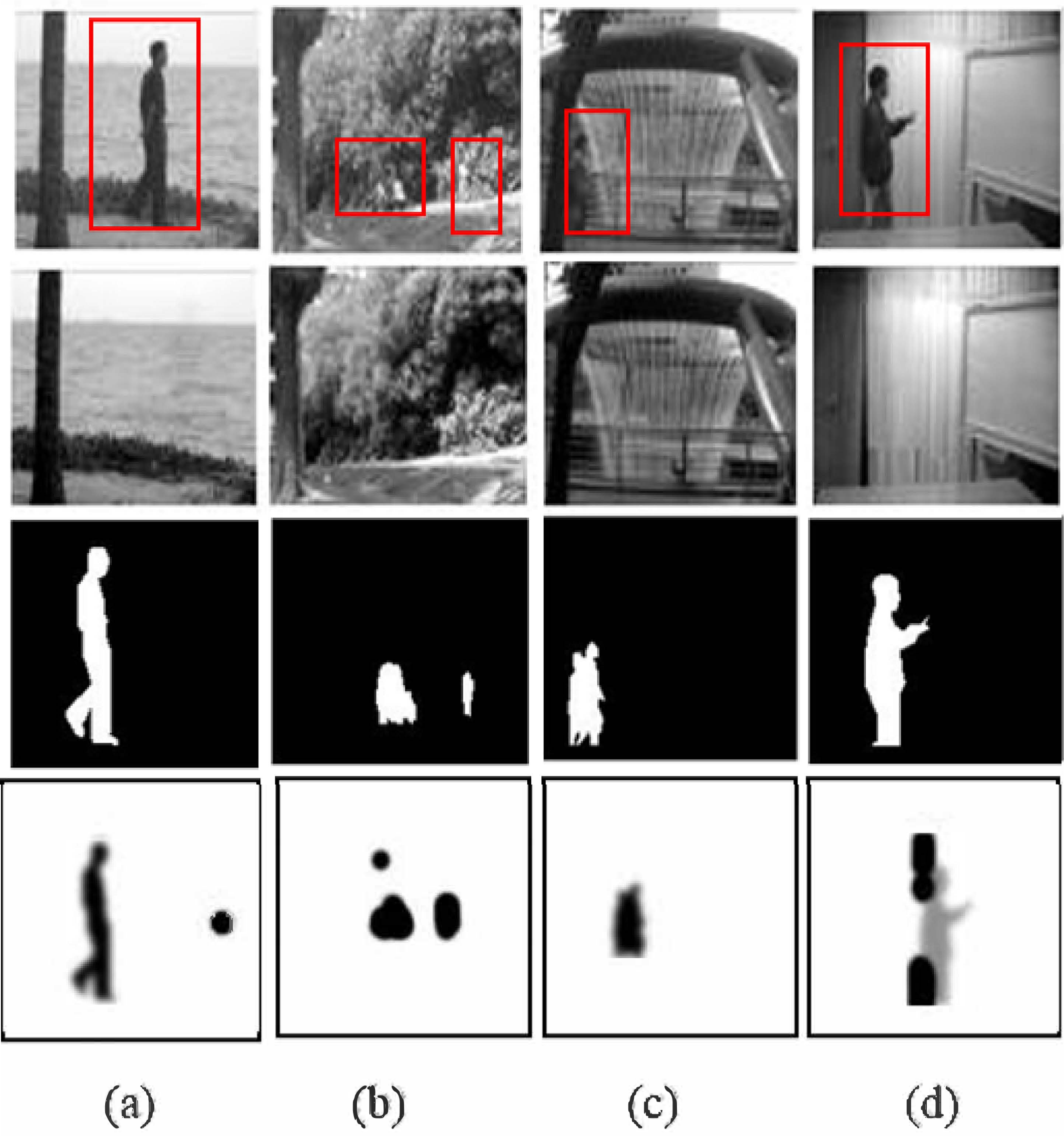}
	\caption{SSIM map of foreground frames:~(a)~{\it Water Surface},~(b)~{\it Waving tree},~(c)~{\it Fountain}, and~(d)~{\it Curtain}.~Top to bottom: Original, background estimated by WLR,~ground truth frame~(size $64\times 80$), SSIM index map~(size $54\times 70$)~for WLR.~The MSSIM are 0.9851, 0.9082, 0.9940,~and~0.9343 respectively.}\label{dyn_for}
\end{figure}
We examine $F$, the foreground  recovered by each method. Since a ground truth mask is available for each video frame, we use a pixel-based measure to form the confusion matrix for the predictive analysis.~In our case, the pixels are represented using 8 bits per sample, and ${\rm M_I}$,~the maximum possible pixel value of the image is 255.~Therefore,a uniform threshold vector linspace(0, ${\rm M}_I$, 100) is used to compare the pixel-wise predictive analysis between each recovered foreground frame and the corresponding ground truth frame.
From the ROC curves in Figure~\ref{roc_curve}, except the {\it Light switch} scenario, the increment in performance of WLR compare to RPCA algorithms does not appear to be substantial.~On the other hand, the qualitative performance of the proposed weighted algorithm in all three scenarios is much superior.~We attribute this to the fact that WLR removes the noise uniformly from the video sequence.~This may lead to an attenuation in performance due to introduced false positives~(see for example Figure~\ref{fp}).


In calculating the PSNR, we perceive the information how the high intensity regions of the image are coming through the noise, and consequently, we pay much less attention to the low intensity regions. This motivated us to remove the noisy components from the recovered foreground, $F$, by using the threshold $\epsilon_1$~(see Algorithm 2), such that we set the components below $\epsilon_1$ in $F$ to 0.~Using this new $F$, we will give the next two quantitative measures. PSNR is calculated using the metric: $10log_{10}\frac{{\rm M}_I^2}{{\rm MSE}}$, where ${\rm MSE}= \frac{1}{mn}\|F(:,i)-G(:,i)\|_2^2$. Conventionally, the higher the PSNR value, the better the reconstruction algorithm.~Figure~\ref{psnr} indicates the PSNR of WLR is superior than the RPCA algorithms.~This can be attributed by the fact that the foreground frames recovered by WLR in all three scenarios, are identical to the ground truth frames.~Hence, they have 0 MSE~(see Figure~\ref{mse}),~resulting in infinity for PSNR. 
Finally we use the mean SSIM~(MSSIM) index to evaluate the overall image quality~\cite{mssim}.~In order to calculate MSSIM of each recovered foreground video frame, we consider a $11\times 11$ Gaussian window with standard deviation~($\sigma$) 1.5.~In Figure~\ref{mssim}, we plot the MSSIM of different methods for all three scenarios.~The MSSIM plot demonstrates that WLR has superior performance over the RPCA algorithms, especially when there is no foreground or static foreground. Moreover, when the video sequence has sensor noise, compression artifacts, and camouflage, for example the {\it Noisy night} scenario, WLR clearly outperforms RPCA.~In Figure~\ref{ssim_map} the SSIM index map of two sample foreground video frames indicate fragmentary foreground recovered by the RPCA algorithms.
\begin{figure}
	\centering  \includegraphics[width=\linewidth]{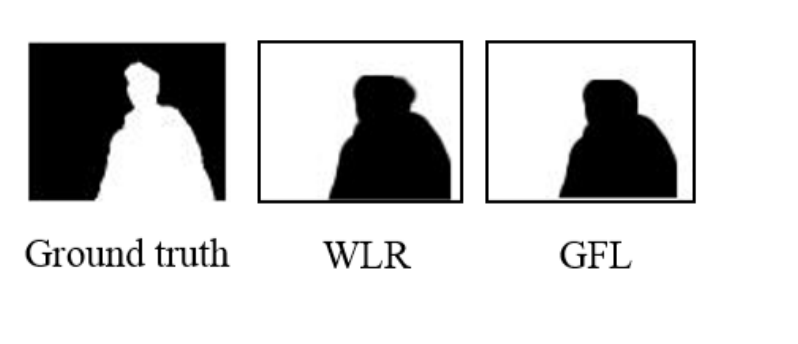}
	\caption{SSIM map of foreground frame of {\it Waving tree} and {\it Basic} scenario. Left to right: Ground truth frame~(size $64\times 80$), SSIM index map~(size $54\times 70$)~for WLR and GFL. MSSIM for WLR and GFL are 0.5018 and 0.5014 respectively.}\label{ssim_map_247}
\end{figure}
\begin{figure}
	\centering  \includegraphics[width=\linewidth, height=1.5in]{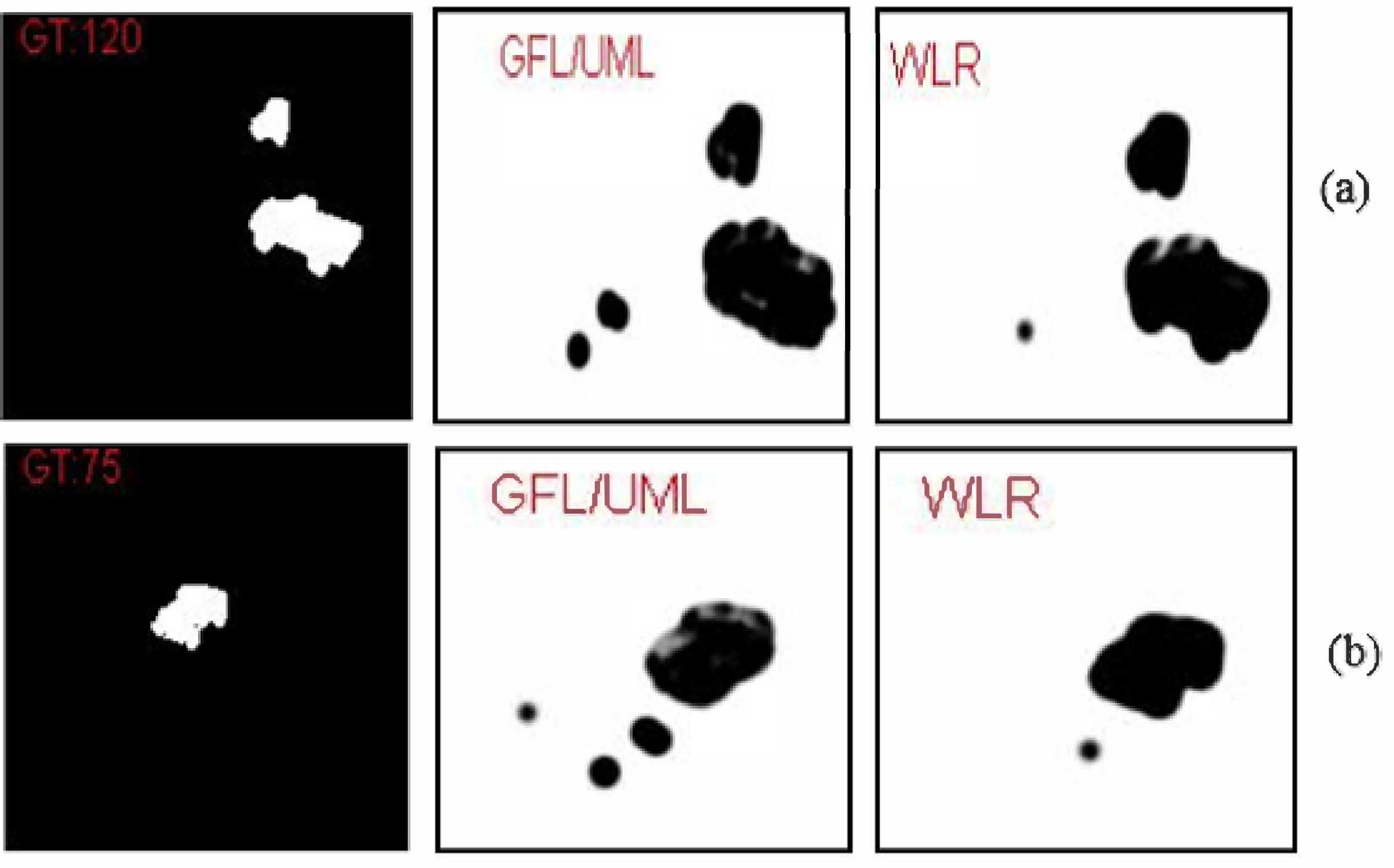}
	\caption{SSIM map of foreground frame of {\it Basic} scenario. Left to right: Ground truth frame~(size $144\times 176$), SSIM index map~(size $134\times 166$)~for WLR and GFL. MSSIM for WLR and GFL are~(a) frame 120:~0.9326 and 0.9244,~(b)~frame 75:~0.9659 and 0.9677 respectively.}\label{ssim_map_st}
\end{figure}
\vspace{-0.1in}
\subsection{Further Experiments on Dynamic Background}
\vspace{-0.05in}
To demonstrate the power of our method on more complex data sets containing dynamic foreground, we perform extensive qualitative and quantitative analysis on the Li data set~\cite{lidata}.~We use four sequences of the data set containing dynamic foreground.~The SSIM index map on all four recovered foreground indicates that WLR performs consistently well on the video sequences containing dynamic background~(see Figure~\ref{dyn_for}). 
\vspace{-0.1in}
\subsection{Comparison with GFL}
\vspace{-0.05in}
We compare the performance of our algorithm with the supervised and upsupervised background subtraction model via Generalized Fused Lasso~of Xin \etal \cite{xin2015}.~Since the choice of $r$ and $k$ are problem specific for our model we have only provided the quantitative comparison on the {\it Waving tree} scene of the Wallflower dataset~\cite{wallflower} and {\it Basic} scenario of the Stuttgart dataset.~Since some frames of the {\it Waving Tree} scenario contain pure background information, Xin \etal used 200 frames as a prior for supervised learning. On the other hand, we used all 286 test frames of the Wallflower sequence to learn the weight and estimate the background without using the exact location of the pure background frames.~From SSIM index map in Figure~\ref{ssim_map_247} it is clear that both methods are very competitive.

~For the {\it Basic} sequence of the Stuttgart dataset we use the unsupervised GFL without using the knowledge of pure background frames.~We use first 200 frames of the {\it Basic} sequence for the unsupervised GFL model and resize the frames as described in the software~\footnote{http://idm.pku.edu.cn/staff/wangyizhou/}~\cite{xin2015}.~For fair comparison we use the same data matrix for WLR.~From SSIM index map in Figure~\ref{ssim_map_st} it is clear that both methods are very competitive with WLR being extraordinarily time efficient than the unsupervised GFL model.~WLR takes approximately 17.75 seconds to conduct the experiment. 

\vspace{-.15in}
\section{Conclusion}
\vspace{-0.08in}
In this paper, we proposed a simple and fast weighted low-rank approximation algorithm for a special family of weights. In addition, we devised an efficient and robust background estimation model and demonstrated its effectiveness on complex video sequences over the existing RPCA algorithms.~The main motivation of the paper is not to propose a background estimation model, rather show how a properly weighted Frobenius norm can be made robust to the outliers, similarly to RPCA and GFL. 

{\small
\bibliographystyle{ieee}
\bibliography{egbib}
}

\end{document}